\documentclass{siamart190516}
\usepackage{lipsum}
\usepackage{amsfonts}
\usepackage{graphicx}
\usepackage{epstopdf}
\usepackage{algorithmic}

\ifpdf
  \DeclareGraphicsExtensions{.eps,.pdf,.png,.jpg}
\else
  \DeclareGraphicsExtensions{.eps}
\fi

\usepackage{xcolor}
\usepackage{color}
\usepackage{tcolorbox}
\usepackage{float}
\usepackage{caption}
\usepackage{amssymb,amsmath,array,graphicx}
\usepackage{graphicx}
\usepackage{textcomp}
\usepackage{enumerate}
\usepackage{epsf}
\usepackage{latexsym,subfigure,paralist}
\usepackage{multirow}
\usepackage{dcolumn}
\newcolumntype{d}[1]{D{.}{.}{#1}}

\newcommand\Tstrut{\rule{0pt}{2.6ex}}         % = `top' strut
\newcommand\Bstrut{\rule[-1.3ex]{0pt}{0pt}}   % = `bottom' strut
\newsiamremark{remark}{Remark}
\newsiamremark{hypothesis}{Hypothesis}
\crefname{hypothesis}{Hypothesis}{Hypotheses}
\newsiamthm{claim}{Claim}

\headers{A New Block Preconditioner for IRK Methods}{M. Rana, V. Howle, K. Long, A. Meek, W. Milestone}

\title{A New Block Preconditioner for \\
Implicit Runge--Kutta Methods \\ 
for Parabolic PDE Problems}
\author{Md Masud Rana\thanks{Corresponding author (co-authors listed alphabetically); 
        Texas Tech University
        (\email{md-masud.rana@ttu.edu})}
        \and
        Victoria E. Howle\thanks{Texas Tech University 
        (\email{victoria.howle@ttu.edu})}
        \and
        Katharine Long\thanks{Texas Tech University
        (\email{katharine.long@ttu.edu})}
        \and
        Ashley Meek\thanks{Adams State University
        (\email{ameek@adams.edu})}
        \and
        William Milestone\thanks{Texas Tech University
        (\email{William.Milestone@ttu.edu})}}
\usepackage{amsopn}
\DeclareMathOperator{\diag}{diag}

\begin{document}
\maketitle
%\maketitle
%\begin{center} \today \end{center}

\begin{abstract}
    A new preconditioner based on a block $LDU$ factorization with algebraic multigrid subsolves for scalability is introduced for the large, structured systems appearing in 
    implicit Runge--Kutta time integration of parabolic partial differential equations. This preconditioner is compared in condition number and eigenvalue distribution, and in numerical experiments with others in the literature: block Jacobi, block Gauss--Seidel, and the optimized
    block Gauss--Seidel method of Staff, Mardal, and Nilssen [{\em Modeling, Identification and Control}, 27 (2006), pp. 109--123].
    %~\cite{staff2006preconditioning}. 
    Experiments are run on two test problems, a $2D$ heat equation and a model advection-diffusion problem, using implicit Runge--Kutta methods with two to seven stages. We find that 
    the new preconditioner outperforms the others, with the improvement becoming more pronounced as spatial discretization is refined and as temporal order is increased. 
\end{abstract}

\section{Introduction}
\label{sec:intro}
Explicit time integrators for parabolic PDEs
are subject to the restrictive timestep limit $h_t \lesssim h_x^2$, so A-stable integrators are essential. It is well known that 
although there are no A-stable explicit linear multistep methods and implicit multistep methods cannot be A-stable beyond order two, there exist A-stable and L-stable implicit Runge--Kutta (IRK) methods at all orders~\cite{hairer1993solving,wanner1996solving}. IRK methods offer an appealing combination of stability and high order; however, these methods are not widely used for PDEs because they lead to large, strongly coupled linear systems. An $s$-stage IRK system has $s$ times as many
degrees of freedom as the systems resulting from backward Euler or implicit trapezoidal rule discretization applied
to the same equation set. Order-optimal preconditioners for such systems have been investigated in a series of papers \cite{mardal2007order,staff2006preconditioning}. In this paper we introduce a new block preconditioner for IRK methods,
based on an LDU factorization. Solves on individual blocks are accomplished using a multigrid algorithm. 
We demonstrate the effectiveness of this preconditioner on two test problems. The first is a simple heat equation, and the second is a model advection-diffusion problem known as the double-glazing problem. 
We find that our preconditioner 
is scalable ($h_x$-independent) and yields better timing results than other preconditioners currently in the literature.

%\begin{itemize}
%\item Our immediate motivation comes wanting to incorporate spatial dynamics
      %into an ecological stoichiometry model.
%\item Ecological stoichiometry is the study of the balance of energy and
      %multiple chemical elements in ecological interactions (Sterner \& Elser,
      %Ecological Stoichiometry, 2002).
%\item For example, predator-prey models of algae-Daphnia systems that
      %incorporate carbon, nitrogen, and phosphorus levels in the environment
      %and effects of light levels.
%\item Incorporating spatial dynamics allows for the fact that light
      %levels are not uniform throughout the lake -- varying with depth.
%\item The researchers on this project are interested in higher order
      %methods and stable methods for stiff systems.
%\item Such methods may also be useful in incompressible flow problems.
%\end{itemize}
%\item IRK allows stable time-stepping with high order
%\item Challenge: Systems arises after discretization are large and coupled, not amenable to Multi Grid (MG) or Incomplete LU (ILU).    

\section{Implicit Runge--Kutta Methods}
\label{sec:IRK}
IRK methods are described in detail in, for example,~\cite{wanner1996solving}.
The general $s$-stage IRK method for $u'=f(t,u)$ requires the solution of $s$ simultaneous equations
$$K_i=f\left(t_n+c_i h_t, u_n+ h_t \sum_{j=1}^s a_{ij}K_j\right)\;\; \mathrm{for}\; i=1\;\mathrm{to}\; s$$
for the stage variables $K_i$, $i=1\;\mathrm{to}\;s$. The stage variables are then used to update the solution $u$,
$$u_{n+1}=u_n+ h_t \sum_{i=1}^s b_i K_i.$$
The coefficients that define a given method are summarized in a Butcher table:
\begin{equation*}
\begin{array}
{c|cccc}
c_1 & a_{11} &\cdots & a_{1s}\\
\vdots & \vdots &\ddots &\vdots\\
c_s &a_{s1} &\cdots & a_{ss} \\
\hline
& b_1 &\cdots & b_s 
\end{array}
\;=\;
\begin{array}
{c|c}
c & A \\
\hline \\[-0.9em]
& b^T.
\end{array}
\end{equation*}

%Dahlquist
%\begin{itemize}
%\item Our goals are higher order and A-stable (or L-stable).
%\item Dahlquist Second Barrier: There are no explicit A-stable linear 
%      multistep methods. \\
%      A-stable implicit linear multistep methods are $2^{nd}$ order at best.
%\item But we have IRK methods that are high order and A or L-stable.
%\end{itemize}

We list some common IRK methods with their order and stability properties in Table \ref{tab:commonIRK}. As our goal is to use higher order L-stable methods, throughout this paper we consider Radau IIA and Lobatto IIIC methods. Of the common L-stable methods, these provide the highest order for a given number of stages, that is, for a given cost. 

\begin{table}[htbp]
\caption{Common IRK Methods}
\centering
\begin{tabular}{l l l}
\hline
IRK Methods & Order & Stability\Tstrut\Tstrut\Bstrut\Bstrut\\
\hline
%Backward Euler & 1 & L-stable\\
%\hline
%Implicit Trapezoidal & 2 & A-stable \\
%\hline
Gauss--Legendre (s) & $2s$ & A-stable\Tstrut\\
%\hline
Radau IA (s) & $2s-1$ & L-stable\\
%\hline
Radau IIA (s) &$2s-1$& L-stable\\
%\hline
Lobatto IIIA (s) & $2s-2$ & A-stable\\
%\hline 
Lobatto IIIC (s) & $2s-2$ & L-stable\\
%\hline 
Miller DIRK (2) & 2 & L-stable\\
%\hline
Miller DIRK (3) & 2 & L-stable\\
%\hline
Crouzeix SDIRK(2) & 2 & L-stable\\
%\hline
Crouzeix SDIRK(3) & 3 & L-stable\Bstrut\\
\hline
\end{tabular}
\label{tab:commonIRK}
\end{table}

\section{Formulation of the Test Problems}
\label{sec:TestProblem}
We consider two test problems. The first is a simple heat equation, and the second is a double-glazing 
advection-diffusion problem. 
\subsection{Heat Equation}
\label{sec:TestProblemHeat}
As a first simple test problem, we consider the heat equation
on some spatial domain $\Omega$ over the time interval $[0,T]$ with
Dirichlet boundary conditions:
\begin{align}\label{eqn:testProblemHeat}
\begin{split}
%u_t &= \nabla ^2 u = f(t,u) \quad \mathrm{in}\; \Omega\times [0,T]\\
u_t &= \nabla ^2 u \quad \mathrm{in}\; \Omega\times [0,T]\\
u &= 0 \quad \mathrm{on}\;\;\partial\Omega\\
u(\bar{x},0) &= u_0(\bar{x}).
\end{split}
\end{align}

We first discretize in time using an IRK method. We write $u_n$ for $u(\bar{x},n h)$ and introduce stage variables $K_i = K_i(\bar{x})$ for $i=1$ to $s$. 
In strong form, the stage equations for an $s$-stage IRK method are:
\begin{align*}
K_i &=\nabla ^2 u_n +  h_t \sum_{j=1}^s a_{ij}\nabla ^2 K_j\quad \mathrm{in}\; \Omega \\
K_i &= 0 \quad \mathrm{on}\;\partial\Omega \qquad\;\; \mathrm{for}\; i=1\;\mathrm{to}\; s.
\end{align*}

We discretize this system of time-independent PDEs using the finite element method. Adopting test functions
$\hat{K}_i\in H_0^1$, we obtain a weak form as usual by multiplying 
the $i$-th stage equation by the test function $\hat{K}_i$, integrating, and applying  
Green's identity. This leads to the equations:
\begin{align*}
\int_{\Omega}\left[\hat{K}_{i}K_{i}+\nabla\hat{K}_{i}\cdot\nabla u_{n}+ h_t\sum_{j=1}^{s}a_{ij}\nabla\hat{K}_{i}\cdot\nabla K_{j}\right]\,d\Omega-\int_{\partial\Omega}\hat{K}_{i}\hat{\mathbf{n}}\cdot\nabla u_{n}\,d\Gamma\\
- h_t\sum_{j=1}^{s}a_{ij}\int_{\partial\Omega}\hat{K}_{i}
\hat{\mathbf{n}}\cdot\nabla K_{j}\,d\Gamma&=0,
\end{align*}
which must hold for all $\hat{K}_i\in H_0^1$.
Note that the boundary terms vanish due to the boundary restriction of 
$\hat{K_i}$ to $H_0^1$.
Therefore we need to find $K_i\in H_0^1$ such that\\
\begin{equation*}
\int_{\Omega}\left[\hat{K}_{i}K_{i}+\nabla\hat{K}_{i}\cdot\nabla u_{n}+ h_t\sum_{j=1}^{s}a_{ij}\nabla\hat{K}_{i}\cdot\nabla K_{j}\right]\,d\Omega = 0 \quad \forall \hat{K_i}\in H_0^1.
\end{equation*}

Converting to weak form and discretizing with finite element 
basis functions $\phi_j$, we get the following linear system, which must be solved for the stage variables $K_i$:
\begin{equation}
\label{block_coeff_mtx}
\begin{bmatrix}
  M+a_{11}  h_t F & a_{12}  h_t F & \cdots & a_{1s}  h_t F \\[0.3em]
  a_{21}  h_t F & M+a_{22}  h_t F & \cdots & a_{2s}  h_t F \\[0.3em]
  \vdots & \vdots &\ddots & \vdots \\[0.3em]
  a_{s1}  h_t F & a_{s2}  h_t F & \cdots & M+a_{ss}  h_t F
\end{bmatrix}
\begin{bmatrix}
{K_1}\\[0.3em]
{K_2}\\[0.3em]
%\bar{K_1}\\[0.3em]
%\bar{K_2}\\[0.3em]
\vdots\\[0.3em]
%\bar{K_S}
{K_s}
\end{bmatrix}
= -
\begin{bmatrix}
%F\bar{u_n}\\[0.3em]
%F\bar{u_n}\\[0.3em]
F{u_n}\\[0.3em]
F{u_n}\\[0.3em]
\vdots\\[0.3em]
%F\bar{u_n}
F{u_n}
\end{bmatrix},
\end{equation}
where the elements of the matrices $M$ and $F$ are given by
\begin{align*}
M_{kl} &= \int \phi_k \phi_l = O\left( h_x^d \right) \\
F_{kl} &= \int \nabla \phi_k \cdot \nabla \phi_l = O \left(h_x^{d-2} \right),
\end{align*}
and where $d$ is the spatial dimension of the problem.
Thus we have an $sN\times sN$ linear system to solve at each time step;
this is the system that we need to precondition.

% [committed by Masud] I wrote this section as a possible development of our preconditioners. 
%
% Notation that I have used:
% $A$ : IRK coefficient matrix
% $\mathcal{A}$ : block matrix from discretization
% $\mathcal{P}$ : general block preconditioners for \mathcal{A}
% $\mathcal{P}_J$ : block Jacobi (diagonal) prec.
% $\mathcal{P}_{GSL}$ : blcok Gauss--Seidel prec.
% $\mathcal{P}_{DU}$ : block LDU based upper triangular prec.
% $\mathcal{P}_{LD}$ : block LDU based lower triangular prec.

\subsection{Double-Glazing Advection-Diffusion Problem}\label{sec:TestProblemAD}
As a second test problem, we consider a model advection-diffusion problem known as the 
double-glazing problem. A time-independent version of this problem is described 
in detail in \cite{elman2005fem}; we modify it here to a time-dependent problem. 
This is a simple model for advective-diffusive transport in a cavity 
$\Omega = [-1,1] \times [-1,1]$, where one wall is hot. The equation and boundary conditions are
\begin{align}\label{eqn:testProblemAD}
\begin{split}
u_t  + \left( \overline{w} \cdot \nabla \right) u 
-\epsilon \nabla^2 u &= 0 \quad \mathrm{in}\; \Omega\times [0,T]\\
u &= 0 \quad \mathrm{on}\;\;\partial\Omega_N \cup \partial\Omega_W  \cup\partial\Omega_S \\
u &= 1 \quad \mathrm{on}\;\;\partial\Omega_E \\
u(\bar{x},0) &= u_0(\bar{x}),
\end{split}
\end{align}
where $\epsilon > 0$ and $ \partial\Omega_N, \; \partial\Omega_W, \; \partial\Omega_S $
and $\partial\Omega_E$ are the North, West, South, and East walls of $\Omega$. For initial conditions we take $u_0(\bar{x})=0$ on $\Omega\backslash \partial\Omega_E$; the solution evolves towards the steady solution described in \cite{elman2005fem}.
The transport is  advection dominated when $\epsilon \ll 1$ and
diffusion dominated when $\epsilon \gtrsim 1$.
For the wind $\overline{w}(\bar{x})$, we use the circulating interior flow described in \cite{elman2005fem}.
%For consistency, we use the same $supg$ stabilization as used in
%\cite{elman2005fem} and the same choices of $\epsilon$ values.
As before, we discretize in time using an IRK method then convert the resulting
time-independent PDEs to weak form and discretize in space using stabilized finite elements. 
This results in a linear system for the stage variables $K_i$ of the same structural form as 
\eqref{block_coeff_mtx}.

\section{Preconditioning IRK methods applied to parabolic PDEs}
\label{sec:Prec}

In this paper, we develop a block upper triangular preconditioner and a block 
lower triangular preconditioner for systems of the form  \eqref{block_coeff_mtx}. 

%\subsection{Development of the Preconditioners}
Let $A$ be the $s$-stage IRK coefficient matrix, that is, the matrix 
comprising the elements in the upper-right block of the Butcher table for 
the given IRK method. 
\begin{equation*}
A = \begin{bmatrix}
a_{11} & a_{12} & \cdots & a_{1s} \\
a_{21} & a_{22} & \cdots & a_{2s} \\
  \vdots & \vdots &\ddots & \vdots \\
  a_{s1} & a_{s2} & \cdots & a_{ss} 
\end{bmatrix}.
\end{equation*}

Given an $m \times n$ matrix $B$ and a matrix $C$, recall that their Kronecker product 
$B \otimes C$ is given by
\begin{equation*}
B\otimes C = \begin{bmatrix}
  b_{11} C & b_{12} C & \cdots & b_{1n} C \\[0.3em]
  b_{21} C & b_{22} C & \cdots & b_{2n} C \\[0.3em]
  \vdots & \vdots &\ddots & \vdots \\[0.3em]
  b_{m1} C & b_{m2} C & \cdots & b_{mn} C
\end{bmatrix}.
\end{equation*}
Using this notation, we can write the block matrix \eqref{block_coeff_mtx} arising from the FEM discretization as
\begin{equation*}
    \mathcal{A} = I_s\otimes M +  h_t\*A\otimes F,
\end{equation*}
where $I_s$ is the $s \times s$ identity matrix.

In this paper, we consider several block preconditioners 
$\mathcal{P}$ 
%for the matrix $\mathcal{A}$ 
of the general form 
\begin{equation}
    \label{eqn:PrecForm}
    \mathcal{P} = I_s\otimes M +  h_t \* \widetilde{A}\otimes F,
\end{equation}
where $\widetilde{A}$ is either a diagonal or triangular matrix referred as the {\em preconditioner coefficient matrix}.
We can, for example, define  
a block Jacobi preconditioner $\mathcal{P}_J$ 
and a block Gauss--Seidel preconditioner $\mathcal{P}_{GSL}$ using this notation as follows:
\begin{equation}
\begin{aligned}
    \mathcal{P}_J &= I_s\otimes M + 
         h_t \widetilde{A}_J \otimes F \\
    \mathcal{P}_{GSL} &= I_s\otimes M + 
         h_t\widetilde{A}_{GSL} \otimes F,
\end{aligned}
\end{equation}
where $\widetilde{A}_J$ is the diagonal part of $A$ 
and $\widetilde{A}_{GSL}$ is the lower triangular part of $A$.
%\[
%\begin{array}{cl}
     %\mathcal{P}_J=I_s\otimes M +  h_t\widetilde{A}_J\otimes F, & \mbox{
     %where } \widetilde{A}_J \mbox{ is the diagonal part of }A\\
    %\mathcal{P}_{GSL}=I_s\otimes M +  h_t\widetilde{A}_{GSL}\otimes F, & \mbox{where } \widetilde{A}_{GSL} \mbox{ is the lower triangular part of }A.
%\end{array}
%\]
%$\mathcal{P}_J=I_s\otimes M +  h_t\widetilde{A}_J\otimes F$, where $\widetilde{A}=\mathrm{diag}(A)$.\\
%And, $\mathcal{P}_{GSL}=I_s\otimes M +  h_t\widetilde{A}_{GSL}\otimes F$, where $\widetilde{A}=\mathrm{tril}(A)$, the lower triangular part of A.
These preconditioners and their order optimality have been extensively studied in Mardal et al.\  \cite{mardal2007order} and in G.A. Staff et al.\ \cite{staff2006preconditioning}. 
% VEH - how do we do references in overleaf?
We now introduce two new $LDU$-based block triangular preconditioners for $\mathcal{A}$.

\subsection{{\em LDU}-based block triangular preconditioners}
It has been shown in \cite{staff2006preconditioning} that if all blocks of $\mathcal{A}$ are well preconditioned by a preconditioner $\mathcal{P}$ of the form \eqref{eqn:PrecForm}, then all of the eigenvalues of the preconditioned system will be clustered, and the condition number of the preconditioned system can be approximated by $\kappa{\left(\mathcal{P}^{-1}\mathcal{A}\right)}\approx \kappa{\left(\widetilde{A}^{-1}A\right)}$,
where $\widetilde{A}$ is the preconditioner coefficient matrix.
In other words, if we make $\widetilde{A}$ an effective preconditioner for the Butcher coefficient matrix $A$, then $\mathcal{P}$ can be expected to be an effective (left) preconditioner for the system $\mathcal{A}$.
(In practice, we use the same motivation in developing left and right preconditioners.)
We begin, therefore, by preconditioning $A$.

For the development of our preconditioners, we assume that the IRK coefficient matrix $A$ is invertible and has nonsingular leading principal submatrices. It is proven in \cite{hairer1993solving} that the Butcher coefficient matrices for Gauss--Legendre, Radau IA, Radau IIA, and Lobatto IIIC are nonsingular for all $s$. 
We have confirmed that the leading principal submatrices are all nonsingular for the coefficient matrices arising from Gauss--Legendre, Lobatto IIIC, Radau IA, and Radau IIA IRK methods for all stages that we use in this paper: $s=2$ through $s=7$. 
Given these conditions, $A$ can be factored into $A=LDU$ without pivoting.
We compute formally the $LDU$ factorization $A = LDU$ without pivoting 
\begin{eqnarray*}
\begin{bmatrix}
a_{11} & a_{12} & \cdots & a_{1s} \\
a_{21} & a_{22} & \cdots & a_{2s} \\
 \vdots & \vdots &\ddots & \vdots \\
a_{s1} & a_{s2} & \cdots & a_{ss} 
\end{bmatrix}
&=& 
\begin{bmatrix}
1 \\
l_{21} & 1 \\
 \vdots & \vdots &\ddots \\
l_{s1} & l_{s2} & \cdots & 1
\end{bmatrix}
\begin{bmatrix}
d_{11} \\
  & d_{22} \\
  & & \ddots \\
  & & & d_{ss}
\end{bmatrix}
\begin{bmatrix}
1 & u_{12} & \cdots & u_{1s} \\
  & 1 & \cdots & u_{2s} \\
  & &\ddots & \vdots \\
  & & & 1
\end{bmatrix} 
\end{eqnarray*}
and consider $\widetilde{A}_{DU} = DU$ and
$\widetilde{A}_{LD} = LD$ as preconditioners for $A$.
Since $A = LDU$, we have 
\begin{align*}
A &= LDU \\
A (DU)^{-1} &= L \\
\text{and} \\
(LD)^{-1} A &= U.
\end{align*}
The eigenvalues of $U$ and $L$ are all one, thus we expect  $\widetilde{A}_{DU}$ to be a good right preconditioner and 
$\widetilde{A}_{LD}$ to be a good left preconditioner for $A$. Since the eigenvalues (though not the condition numbers) of the product of two matrices do not depend on the order of the product, we consider both matrices as left and right preconditioners in our numerical experiments. 

Using the condition number motivation from \cite{staff2006preconditioning}, we now 
define block upper and lower triangular preconditioners for $\mathcal{A}$ of the form:
\begin{equation}
\begin{aligned}
\label{eqn:PLDandPDU}
\mathcal{P}_{DU} &= 
    I_s\otimes M  +  h_t\widetilde{A}_{DU}\otimes F \\
\mathcal{P}_{LD} &= 
    I_s\otimes M +  h_t\widetilde{A}_{LD}\otimes F.
\end{aligned}
\end{equation}

\subsection{Application of the Preconditioners}
Applying the block Jacobi preconditioner $\mathcal{P}_J$ to a vector $v$ is straightforward, and for an $s$-stage IRK scheme involves $s$ subsolves. The subsolve for the $j$th block has the form 
\begin{equation*}
    \left( M + h_t \widetilde{A}_{jj} F \right) w_j = 
    v_j.
\end{equation*}
In practice and for efficiency, we use a single algebraic multigrid (AMG) V-cycle for each subsolve for the block $w_j$.
% VEH - is it clear enough that v and w are blocked similarly and v_j and w_j are the jth blocks of those vectors?

Application of the block upper triangular preconditioner 
$\mathcal{P}_{DU}$ is done via back substitution and again involves $s$ subsolves. The subsolve on the $j$th block in this case is of the form 
\begin{equation*}
    \left( M + h_t \widetilde{A}_{jj} F \right) w_j
       = \left[ v_j - \sum_{k = j+1}^s h_t \widetilde{A}_{jk} F w_k \right].
\end{equation*}
In practice, we again use a single AMG V-cycle for each subsolve.

Application of the block lower triangular preconditioners 
$\mathcal{P}_{GSL}$ and  $\mathcal{P}_{LD}$ is done via forward substitution and again involves $s$ subsolves. 
The subsolves on the $j$th block in this case are of the form 
\begin{equation*}
    \left( M + h_t \widetilde{A}_{jj} F \right) w_j
       = \left[ v_j - \sum_{k = 1}^{j-1} h_t \widetilde{A}_{jk} F w_k \right].
\end{equation*}
As with the previous preconditioners, we use a single AMG V-cycle for each subsolve.

\subsection{Analysis of the Preconditioners}
The preconditioners have been chosen because of their intended effect on the 
condition number of the preconditioned system. Although, as expected, the condition number 
is not always a good indicator of performance for GMRES. 
In Table \ref{table:CondNumHeat}, 
we examine the 2-norm condition number for the unpreconditioned system $\mathcal{A}$ 
and for $\mathcal{A}$
preconditioned on the right by  $\mathcal{P}_J$,
$\mathcal{P}_{GSL}$, $\mathcal{P}_{DU}$, and
$\mathcal{P}_{LD}$.
The system $\mathcal{A}$ is from a 2D heat equation using an $s$-stage Radau IIA IRK method for $s$ ranging from 2 to 7. 
In all of these results, we have chosen $h_x = 2^{-3}$ and $h_t= h_x^{\frac{p+1}{2s-1}}$, where $p=2$ is the degree of the Lagrange polynomial basis functions in space. The preconditioners are constructed exactly.
All of the preconditioners significantly reduce the condition number, with $\mathcal{P}_{GSL}$ and 
$\mathcal{P}_{LD}$ giving the most improvement. 
$\mathcal{P}_{GSL}$ gives the greatest improvement on the condition number for the lower-order stages $s=2$ and $s=3$, and $\mathcal{P}_{LD}$ gives the greatest improvement for the higher-order stages $s=3$ through $s=7$.
% first cond number table
\begin{table}[htbp]
\begin{center}
\caption{Condition numbers of right-preconditioned matrices with preconditioners $\mathcal{P}_J^{-1}$, $\mathcal{P}_{GSL}^{-1}$, $\mathcal{P}_{DU}^{-1}$, and $\mathcal{P}_{LD}^{-1}$ applied to a 2D heat equation with $s$-stage Radau IIA methods. Here, $h_x=2^{-3}$ and $ h_t= h_x^{\frac{p+1}{2s-1}}$, where $p=2$ is the degree of the Lagrange polynomial basis functions in space. Preconditioners are constructed exactly.}
%\footnotesize

\begin{tabular}{c c c c c c}
%\begin{tabular}{c d{4.2} d{2.2} d{2.2} d{2.2} d{2.2}}
	\hline
    %   &         & GMRES    & GMRES      & GMRES        & GMRES \\
    s  & $\kappa{(\mathcal{A})}$       & $\kappa{(\mathcal{A}\mathcal{P}_J^{-1})}$ &$\kappa{(\mathcal{A}\mathcal{P}_{GSL}^{-1})}$ & $\kappa{(\mathcal{A}\mathcal{P}_{DU}^{-1})}$ & $\kappa{(\mathcal{A}\mathcal{P}_{LD}^{-1})}$\Tstrut\Tstrut\Tstrut\Tstrut\Bstrut\Bstrut\Bstrut\Bstrut  \\
		\hline
		2   & 240.37     & 3.23 & 1.75  & 5.32    & 2.48\Tstrut  \\
		%\hline
		3  & 502.53    & 5.66  & 2.58  & 11.18    & 2.66 \\
		%\hline
		4  & 746.23   & 8.54 & 3.63  & 18.23   & 3.04 \\
		%\hline
		5  & 959.16   & 11.76 & 5.08  & 26.53   & 3.21\\
		6  & 1137.24   & 15.23 & 7.13  & 35.97   & 3.50\\
		7  & 1281.47   & 18.90 & 10.05  & 46.48   & 3.67\Bstrut\\
    \hline
\end{tabular}
\label{table:CondNumHeat}
\end{center}
\end{table}
	
In Table \ref{table:CondNumDG}, 
we examine the 2-norm condition number for the unpreconditioned system $\mathcal{A}$ arising from the double-glazing problem, 
for $\mathcal{A}$ preconditioned on the left by
$\mathcal{P}_{GSL}$ and $\mathcal{P}_{LD}$, and for $\mathcal{A}$ preconditioned on the right by  
$\mathcal{P}_{GSL}$ and $\mathcal{P}_{LD}$.
We have set $\epsilon = 0.005$, which corresponds to a weakly advection-dominated problem. 
We are using an $s$-stage Radau IIA IRK method for $s$ ranging from 2 to 7. 
In all of these results, we have chosen $h_x = 2^{-4}$ and $h_t= h_x^{\frac{p+1}{2s-1}}$, where $p=1$ is the degree of the Lagrange polynomial basis functions in space. The preconditioners are constructed exactly.
As with the heat problem, both of the preconditioners significantly reduce the condition number. 
In this table, we omit results for the block Jacobi preconditioner and our upper triangular preconditioner since, as in the heat test problem, these did not perform as well as $\mathcal{P}_{GSL}$ and 
$\mathcal{P}_{LD}$. 
The $\mathcal{P}_{LD}$ preconditioner 
applied on the right gives greater improvement than 
$\mathcal{P}_{GSL}$ applied on the right for the higher-order stages $s=4$ through $s=7$.
$\mathcal{P}_{LD}$  applied on the left gives the greatest improvement in condition number for all stages. 

\begin{table}[htbp]
\begin{center}
\caption{Condition numbers of left-preconditioned and right-preconditioned matrices with preconditioners 
$\mathcal{P}_{GSL}^{-1}$ and $\mathcal{P}_{LD}^{-1}$ applied to a 2D double-glazing problem with $\epsilon = 0.005$, and
with $s$-stage Radau IIA methods. Here, $h_x=2^{-4}$ and $ h_t= h_x^{\frac{p+1}{2s-1}}$, 
where $p=1$ is the degree of the Lagrange polynomial basis functions in space. 
Preconditioners are constructed exactly.}
%\footnotesize
\begin{tabular}{c c c c c c c}
	\hline
    %   &         & GMRES    & GMRES      & GMRES        & GMRES \\
    \multicolumn{2}{c}{} & \multicolumn{2}{c}{$\kappa{(\mathcal{P}^{-1}\mathcal{A})}$} & & \multicolumn{2}{c}{$\kappa{(\mathcal{A}\mathcal{P}^{-1})}$}\Tstrut\Bstrut\\
    \cline{3-4}\cline{6-7}
    s  & $\kappa{(\mathcal{A})}$ &$\kappa{(\mathcal{P}_{GSL}^{-1}\mathcal{A})}$ &$\kappa{(\mathcal{P}_{LD}^{-1}\mathcal{A})}$   & &$\kappa{(\mathcal{A}\mathcal{P}_{GSL}^{-1})}$ 
        & $\kappa{(\mathcal{A}\mathcal{P}_{LD}^{-1})}$\Tstrut\Tstrut\Tstrut\Tstrut\Bstrut\Bstrut\Bstrut\Bstrut \\
		\hline
		%2   & 631.32  & 1.67 & 1.26 & & 1.73 & 2.40  \Tstrut\\
		2   & 631.32  & 1.67 & 1.26 & & 1.73 & 2.40  \\
		3   & 1237.60  &2.67  & 1.53 & & 2.51 & 2.57 \\
		4   & 1802.02  &4.14  & 1.80 & & 3.54 & 2.92 \\
		5   & 2310.65  &6.47  & 2.05 & & 4.92 & 3.07 \\
		6   & 2762.17  & 10.23  & 2.25 & & 6.83 & 3.32 \\
		%7   & 3155.26  & 16.30 & 2.44 & & 9.55 & 3.46 \Bstrut\\
		7   & 3155.26  & 16.30 & 2.44 & & 9.55 & 3.46 \\
		%\hline
    \hline
\end{tabular}
\label{table:CondNumDG}
\end{center}
\end{table}

\begin{figure}[htb]
\begin{center}
\subfigure[]{\includegraphics[width=.45\linewidth]{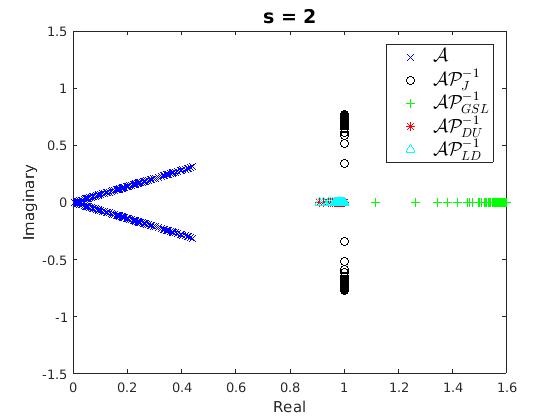}}
\subfigure[]{\includegraphics[width=.45\linewidth]{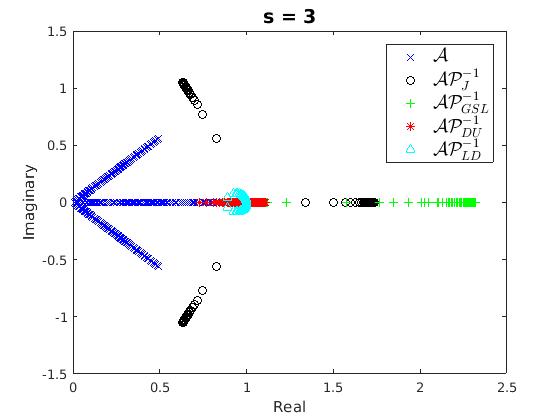}}
\subfigure[]{\includegraphics[width=.45\linewidth]{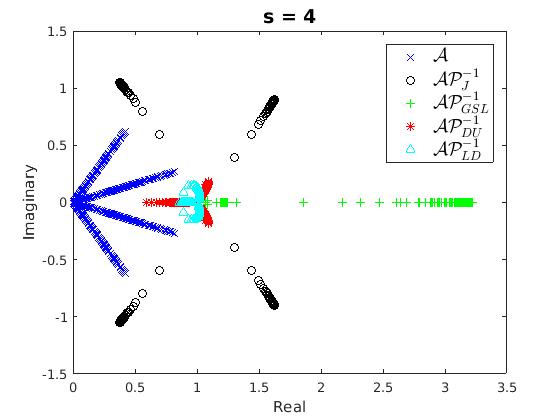}}
\subfigure[]{\includegraphics[width=.45\linewidth]{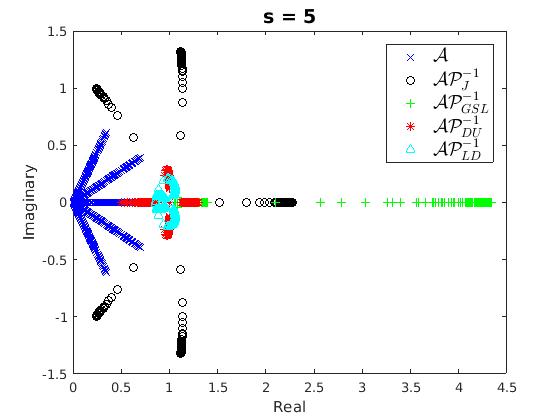}}
\caption{Eigenvalues of the matrices $\mathcal{A}$, $\mathcal{A}\mathcal{P}_J^{-1}$, $\mathcal{A}\mathcal{P}_{GSL}^{-1}$, $\mathcal{A}\mathcal{P}_{DU}^{-1}$, and $\mathcal{A}\mathcal{P}_{LD}^{-1}$ for 2D heat problem with Radau IIA $s=2$ (a), $s=3$ (b), $s=4$ (c), and $s=5$ (d). The $x$-axis is the real part and the $y$-axis is the imaginary part of the eigenvalue. Preconditioners are constructed exactly.}

\label{fig:PrecEwsHeat}
\end{center}
\end{figure}

\begin{figure}[htb]
\begin{center}
\subfigure[]{\includegraphics[width=.45\linewidth]{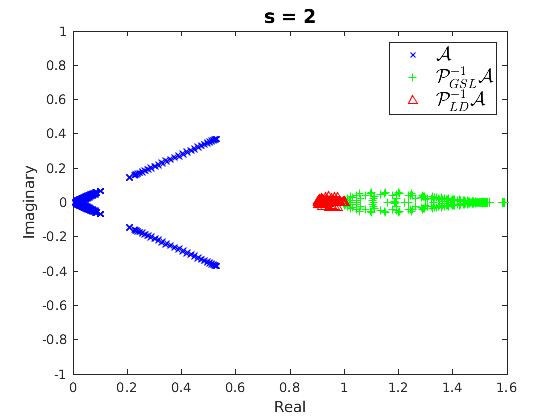}}
\subfigure[]{\includegraphics[width=.45\linewidth]{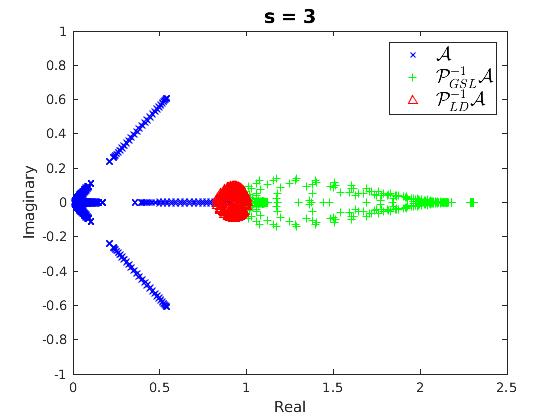}}
\subfigure[]{\includegraphics[width=.45\linewidth]{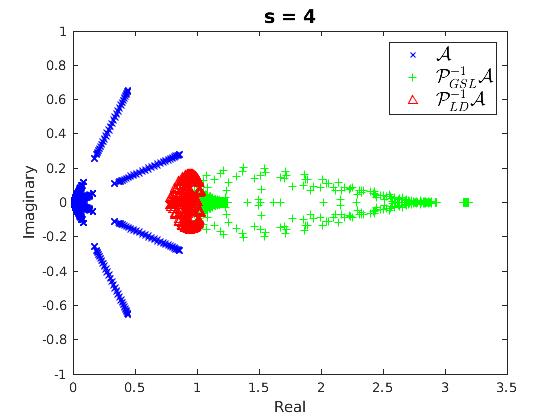}}
\subfigure[]{\includegraphics[width=.45\linewidth]{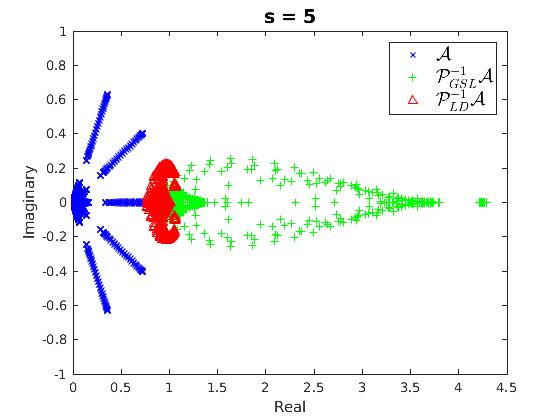}}
\caption{Eigenvalues of the matrices $\mathcal{A}$,  $\mathcal{P}_{GSL}^{-1} \mathcal{A}$, and $\mathcal{P}_{LD}^{-1} \mathcal{A}$ for 2D double-glazing problem with $\epsilon = 0.005$, and with Radau IIA $s=2$ (a), $s=3$ (b), $s=4$ (c), and $s=5$ (d). The $x$-axis is the real part and the $y$-axis is the imaginary part of the eigenvalue. Preconditioners are constructed exactly.}
\label{fig:PrecEwsDG}
\end{center}
\end{figure}

With GMRES, the condition number does not always tell the whole story, and it is also important to consider how well the preconditioner clusters the eigenvalues. 
In Figures \ref{fig:PrecEwsHeat} and \ref{fig:PrecEwsDG}, we show eigenvalues plotted 
in the complex plane for the matrix $\mathcal{A}$ and preconditioned systems. 
In Figure \ref{fig:PrecEwsHeat} the matrix $\mathcal{A}$ is from the 2D heat problem 
and this system is right preconditioned with 
$\mathcal{P}_J$,
$\mathcal{P}_{GSL}$, $\mathcal{P}_{DU}$, and
$\mathcal{P}_{LD}$.
In Figure \ref{fig:PrecEwsDG}, $\mathcal{A}$ is from the 
2D double-glazing problem with $\epsilon = 0.005$ and this system is left preconditioned with 
$\mathcal{P}_{GSL}$ and $\mathcal{P}_{LD}$.
In both figures, we are using 
an $s$-stage Radau IIA IRK method for $s$ ranging from 2 to 5, and 
we have again chosen $h_x = 2^{-3}$ and $h_t= h_x^{\frac{p+1}{2s-1}}$, 
where $p$ is the degree of the Lagrange polynomial basis functions in space; 
$p=2$ for the heat problem and $p=1$ for the double-glazing problem. 
The preconditioners are constructed exactly.
In both cases, the original matrix $\mathcal{A}$ has a number of eigenvalues very near zero. 
All of the preconditioners succeed in clustering the eigenvalues farther away from zero.  
In the heat problem with $s=2$, $\mathcal{P}_{GSL}$ and $\mathcal{P}_{LD}$ give the tightest 
clustering and have very similar eigenvalues.  But as $s$ increases, 
$\mathcal{P}_{LD}$ maintains tighter clustering than
$\mathcal{P}_{GSL}$. This fact is reflected in our numerical results presented in the next section.
In the double-glazing problem, the original matrix $\mathcal{A}$ also has a number of eigenvalues
very near zero. In this problem,  both preconditioners cluster the eigenvalues farther from zero. Our $\mathcal{P}_{LD}$ preconditioner gives the tightest clustering for all stages $s$.

% first ew figure

% DG ew figure

\subsection{Comparison between $\mathcal{P}_{GSL}$ with optimal coefficients and $\mathcal{P}_{LD}$ preconditioners}
In G.A. Staff et al.\ \cite{staff2006preconditioning}, the authors develop an optimal block lower triangular 
Gauss--Seidel preconditioner.
This block Gauss--Seidel preconditioner is derived such that the preconditioner coefficients 
$\widetilde{A}_{GSL}$ are computed from the following optimization problem
 %In this section we compare our $LDU$ based block lower triangular preconditioner $P_{LD}$ with the optimal block lower triangular Gauss--Seidel preconditioner introduced in G.A. Staff et al.\ \cite{staff2006preconditioning}. Let $\widetilde{\mathcal{P}}_{GSL}$ be a GSL preconditioner such that the preconditioner coefficients $\widetilde{A}_{GSL}$ is computed from the following optimization problem
\begin{equation}
\begin{array}{ll}
    \displaystyle{ \min_{\widetilde{A}_{GSL}} } &\kappa(\widetilde{A}_{GSL}^{-1}A)\Bstrut\Bstrut\\
    \mathrm{s.t.}  &\diag{(\widetilde{A}_{GSL})}=\diag{(A)}.
\end{array}
\label{eqn:GSLopt}
\end{equation}	
That is, their optimal preconditioner, which we will refer to as $\widetilde{\mathcal{P}}_{GSL}$, uses coefficients that optimize the condition number of the (left) preconditioned system for a block lower triangular matrix subject to the given constraint. 

In Table \ref{tab:optimizedGSL}, 
we examine the 2-norm condition number for 
$\mathcal{A}$ left preconditioned by 
$\widetilde{\mathcal{P}}_{GSL}$ and $\mathcal{P}_{LD}$ 
for 1D and 2D heat problems. 
For $\widetilde{\mathcal{P}}_{GSL}$, we use the 
optimized $\widetilde{A}_{GSL}$ coefficients given in 
\cite{staff2006preconditioning}.
We fix $h_x=2^{-3}$ and 
$h_t= h_x^{\frac{p+1}{q}}$, where $p=2$ is the degree of the Lagrange polynomial 
basis functions in space and $q$ is the order of the corresponding IRK method 
(that is, $q = 2s-1$ for Radau IIA and $q = 2s -2$ for Lobatto IIIC).
We use Radau IIA with $s=2$ through $s=6$ and Lobatto IIIC with $s=2$ through $s=4$.
For the Lobatto IIIC problems, we see that 
$\widetilde{\mathcal{P}}_{GSL}$ consistently yields a lower condition number than $\mathcal{P}_{LD}$ 
for all stages $s$.
For the Radau IIA problems, however, the $\mathcal{P}_{LD}$ preconditioner yields a lower condition number 
than the optimized preconditioner  $\widetilde{\mathcal{P}}_{GSL}$
for all stages $s$. 
We note that in \cite{staff2006preconditioning}, the authors comment that 
they may not obtain the global minimum to \eqref{eqn:GSLopt} with the optimization process 
they use. The values used in $\mathcal{P}_{LD}$ give a lower condition number.
%use a Nelder-Mead  algorithm initialized with values from $A_{GSL}$ and 
% second heat cond number table
\begin{table}[htbp]
\begin{center}
\caption{Condition numbers of the left-preconditioned system with preconditioners $\widetilde{\mathcal{P}}_{GSL}$ ($\mathcal{P}_{GSL}$ with optimal coefficients) and $\mathcal{P}_{LD}$ for various IRK methods applied to $1D$ and $2D$ heat problems. Here, $h_x=2^{-3}$ and $ h_t= h_x^{\frac{p+1}{q}}$, where $p=2$ is the degree of the Lagrange polynomial in space
and $q$ is the order of the corresponding IRK method
(that is, $q = 2s-1$ for Radau IIA and $q = 2s -2$ for Lobatto IIIC).
Preconditioners are constructed exactly.}

\begin{tabular}{c c c c c c}
\hline
%\\
\multicolumn{3}{c}{1D} & & \multicolumn{2}{c}{2D}\Tstrut\Tstrut\Tstrut\\
\cline{2-3}\cline{5-6}\\[-1.0em]
& $\kappa{(\widetilde{\mathcal{P}}_{GSL}^{-1}\mathcal{A})}$ & $\kappa{(\mathcal{P}_{LD}^{-1}\mathcal{A})}$ &\,& $\kappa{(\widetilde{\mathcal{P}}_{GSL}^{-1}\mathcal{A})}$ & $\kappa{(\mathcal{P}_{LD}^{-1}\mathcal{A})}$\Tstrut\Tstrut\Tstrut\Bstrut\Bstrut\Bstrut\\
\hline
R IIA 2 &1.58 & 1.26 &\,& 1.60   & 1.26\Tstrut \\
R IIA 3 & 1.94 & 1.50 &\,& 1.93  & 1.52 \\
R IIA 4 & 2.18 & 1.74 &\,& 2.17 & 1.77 \\
R IIA 5 & 2.37 & 1.95 &\,& 2.40 & 1.98 \\
R IIA 6 & 3.01 & 2.14 &\,& 3.18 & 2.18\Bstrut\\
\hline
L IIIC 2 &1.48 &2.68 &\,& 1.54 &2.78\Tstrut \\
L IIIC 3 &3.20 &6.68 &\,& 3.32 &6.94 \\
L IIIC 4 &4.92 &10.43 &\,& 5.24 &10.93\Bstrut \\
\hline
\end{tabular}
\label{tab:optimizedGSL}
\end{center}
\end{table}

%G 2 & 1.33 & 1.19 &\, &1.34 & 1.20\Tstrut \\
%G 3 & 1.55 &1.36 &\,& 1.54 & 1.37 \\
%G 4 & 1.63 &1.53 &\,& 1.67 & 1.54 \\
%G 5 & 1.91 & 1.68 &\,& 2.01 & 1.70\Bstrut \\
%\hline

\begin{figure}[htbp]
\begin{center}
\subfigure[]{\includegraphics[width=.45\linewidth]{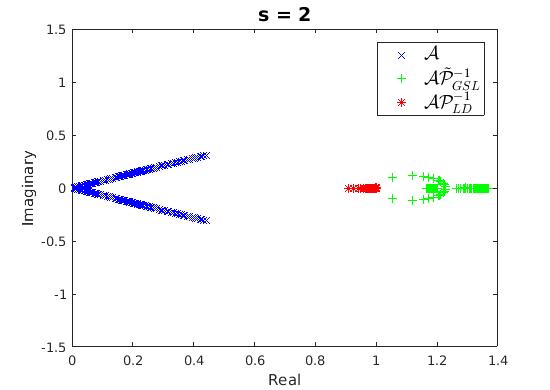}}
\subfigure[]{\includegraphics[width=.45\linewidth]{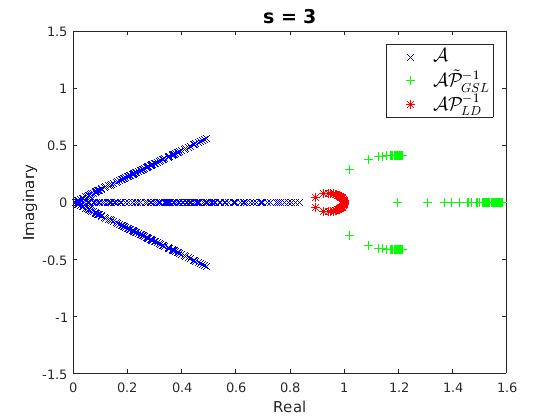}}
\subfigure[]{\includegraphics[width=.45\linewidth]{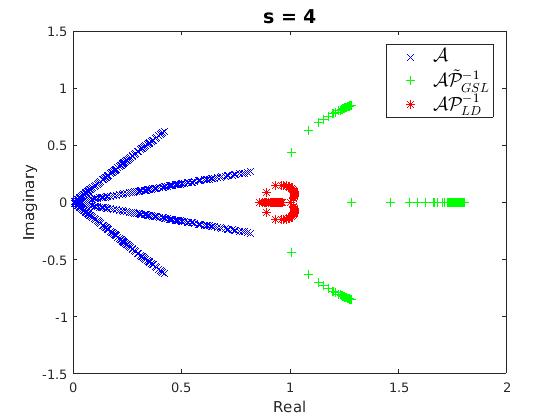}}
\subfigure[]{\includegraphics[width=.45\linewidth]{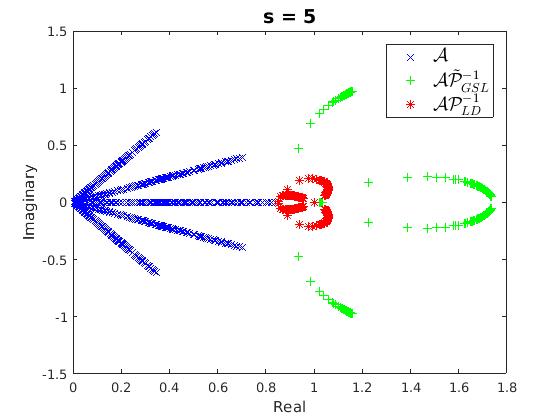}}
\caption{Eigenvalues of the matrix $\mathcal{A}$, $\mathcal{A}\widetilde{\mathcal{P}}_{GSL}^{-1}$, and $\mathcal{A}\mathcal{P}_{LD}^{-1}$ for 2D heat problem with Radau IIA $s=2$ (a), $s=3$ (b), $s=4$ (c), and $s=5$ (d). The $x$-axis is the real part and the $y$-axis is the imaginary part of the eigenvalue. Preconditioners are constructed exactly.}
\label{fig:ewsOptGSLRadau}
\end{center}
\end{figure}

\begin{figure}[htbp]
\begin{center}
\subfigure[]{\includegraphics[width=.45\linewidth]{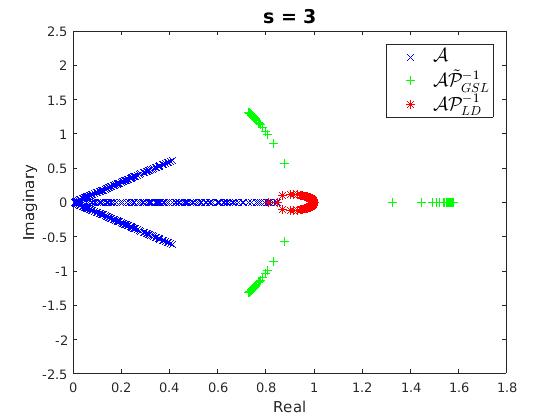}}
\subfigure[]{\includegraphics[width=.45\linewidth]{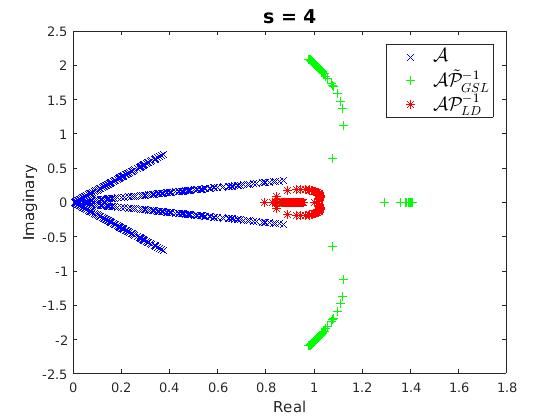}}
\caption{Eigenvalues of the matrix $\mathcal{A}$, $\mathcal{A}\widetilde{\mathcal{P}}_{GSL}^{-1}$, and $\mathcal{A}\mathcal{P}_{LD}^{-1}$ for 2D heat problem with Lobatto IIIC stages $s=3$ (a) and $s=4$ (b). The $x$-axis is the real part and the $y$-axis is the imaginary part of the eigenvalue. Preconditioners are constructed exactly.}
\label{fig:ewsOptGSLLobatto}
\end{center}
\end{figure}

Again we note that for GMRES, a lower condition number does not always indicate a superior preconditioner. 
In Figures \ref{fig:ewsOptGSLRadau} and 
\ref{fig:ewsOptGSLLobatto} we see that $\mathcal{P}_{LD}$ more effectively clusters the eigenvalues for both Radau IIA and Lobatto IIIC for stages $s=2,3,4$ and $5$.

% second ew figure for heat
	
% 3rd ew figure	for heat

\section{Numerical Experiments}
\label{sec:Numerical}
To test the performance of our preconditioners we consider the two model problems described in Section \ref{sec:TestProblem}: a $2D$ heat equation problem and the $2D$ advection-diffusion double-glazing problem. For both we use 
$s$-stage Radau IIA or Lobatto IIIC methods for time discretization.
For the heat equation 
we discretize in space using Galerkin finite elements with piecewise quadratic basis functions on triangular meshes. For the double-glazing problem we use streamwise upwind Petrov-Galerkin discretization~\cite{elman2005fem} with piecewise linear basis functions, again on triangular meshes. Dirichlet boundary conditions are imposed using 
Nitsche's method~\cite{stenberg1995some,juntunen2009nitsche}.

Production of the finite element matrices is carried out with the Sundance finite element toolkit~\cite{long_unified_2010}. The method of manufactured solutions is used to produce an exact solution to each linear subproblem so that we can check the error in addition to the residual. 

A common practice in analysis of methods for spatiotemporal problems is to choose $h_t\sim h_x$ and increase the polynomial
order of spatial discretization in lockstep with the order of time discretization. Here, we use a different approach. We
use low-order elements ($p=1$ or $p=2$) throughout, and then for each $h_x$ we choose a timestep $h_t$ such that the spatial discretization error and temporal global truncation error are comparable. With a Radau $s$-stage method, we set $h_x^{p+1}= h_t^{2s-1}$ so that the timestep is $h_t= h_x^{\frac{p+1}{2s-1}}$, and 
with a Lobatto IIIC method, we set $h_x^{p+1}= h_t^{2s-2}$ so that the timestep is $h_t= h_x^{\frac{p+1}{2s-2}}$.
By holding the spatial order $p$ fixed, we are simulating conditions in which a modeler uses a high order integrator to obtain results of a specified accuracy with larger timestep and smaller computational cost. 

We solve the linear system \eqref{block_coeff_mtx} using preconditioned GMRES. We construct our preconditioners as described in Section \ref{sec:Prec} using one AMG V-cycle for each subsolve. We use AMG from the IFISS software package \cite{ifiss}. When applied to the advection-diffusion equation, multilevel algorithms and smoothers must be handled with care~\cite{elman2005fem}. In these experiments we use the strategy of ensuring that the SUPG discretization is stabilized on the coarse mesh, and use the smoothed aggregation smoother from IFISS. All results are computed in MATLAB on a machine with an Intel Core i7 1.80 GHz (Turbo Boost up to 4.00 GHz) processor and 8.00 GB RAM.

In Section \ref{sec:Prec}, we examined condition numbers and eigenvalue
distributions for the preconditioners 
$\mathcal{P}_{J}$, 
$\mathcal{P}_{GSL}$,
$\widetilde{\mathcal{P}}_{GSL}$, 
$\mathcal{P}_{DU}$,
and $\mathcal{P}_{LD}$
applied to the $2D$ heat equation \eqref{eqn:testProblemHeat} and the
double-glazing advection-diffusion problem \eqref{eqn:testProblemAD}.
In this section, we investigate the convergence performance of our preconditioner compared to the others. 
We start with the $2D$ heat equation \eqref{eqn:testProblemHeat} 
in Section \ref{sec:NumericalHeat}
and compare iteration counts and timing results for GMRES preconditioned on the right by 
$\mathcal{P}_{J}$, 
$\mathcal{P}_{GSL}$,
$\mathcal{P}_{DU}$, and
$\mathcal{P}_{LD}$
for various problem sizes and varying orders of IRK methods.
We then examine the robustness of the 
$\mathcal{P}_{GSL}$ and  $\mathcal{P}_{LD}$ preconditioners for a fixed 
spatial resolution as we vary the time step size. 
We also compare our preconditioner with the optimized preconditioner 
$\widetilde{\mathcal{P}}_{GSL}$.  We conclude the heat equation results
by examining the relative errors achieved with each of 
the preconditioners. 
In Section \ref{sec:NumericalAD} we give iteration count and 
timing results for 
the best performing of the preconditioners, 
$\mathcal{P}_{GSL}$ and $\mathcal{P}_{LD}$,
applied to the double-glazing advection-diffusion problem
\eqref{eqn:testProblemAD}. As with the heat equation results, we first 
compare methods for varying spatial resolutions and various order IRK
methods. We then fix the spatial resolution and compare the robustness 
the methods as we vary the time step size. 

%Preview of numerical experiments: \textcolor{red}{flesh out this list}
%\begin{itemize}
    %\item Eigenvalues and condition numbers (\textcolor{blue}{should the results above be moved here?})
    %\item Convergence performance
    %\item Varying timestep
%\end{itemize}

\subsection{Heat Equation Results}
\label{sec:NumericalHeat}
We first test the performance of the preconditioners 
$\mathcal{P}_J$, $\mathcal{P}_{GSL}$, $\mathcal{P}_{DU}$, and $\mathcal{P}_{LD}$, on our $2D$ heat problem \eqref{eqn:testProblemHeat}.

\subsubsection{Comparison results between $\mathcal{P}_J$, $\mathcal{P}_{GSL}$, $\mathcal{P}_{DU}$, and $\mathcal{P}_{LD}$ preconditioners}

In Table \ref{tab:manufacturedResultsRadauIIA} we report iteration counts and timing results for right-preconditioned GMRES to converge with relative tolerance $1.0 \times 10^{-8}$ for varying mesh sizes $h_x$ and varying number of stages (varying order).
Throughout, all of the triangular preconditioners outperform the block Jacobi preconditioner $\mathcal{P}_J$, as expected. The $LDU$-based block upper triangular preconditioner $\mathcal{P}_{DU}$ and the block Gauss--Seidel preconditioner $\mathcal{P}_{GSL}$ perform similarly, with $\mathcal{P}_{GSL}$ performing slightly better, especially for the higher stage methods. 
However, overall, the $LDU$-based block lower triangular preconditioner $\mathcal{P}_{LD}$ performs better than all of the other preconditioners in both iteration count and timing. 
For 2-stage Radau IIA, $\mathcal{P}_{LD}$ and $\mathcal{P}_{GSL}$ have similar performance, with 
$\mathcal{P}_{LD}$ performing slightly better on the smaller problems and  $\mathcal{P}_{GSL}$ performing slightly better on the larger problems. 
But as the number of stages increases (that is, with increasing order), the performance of 
$\mathcal{P}_{LD}$ improves over the other preconditioners giving consistently lower iteration counts and timing.
For a given number of stages, all of the triangular preconditioners scale well with problem size, exhibiting little to no growth in iteration count with increasing problem size. However the $\mathcal{P}_{LD}$ preconditioner has the least increase in cost as the number of stages increases.

\begin{table}[htbp]
	\begin{center}
	 \caption{Iteration counts and elapsed time (times in seconds are shown in parentheses) for right-preconditioned GMRES to converge with relative residual tolerance $1.0 \times 10^{-8}$ for a $2D$ heat problem with $s$-stage Radau IIA methods with preconditioners $\mathcal{P}_J$, $\mathcal{P}_{GSL}$, $\mathcal{P}_{DU}$, and $\mathcal{P}_{LD}$. Here we choose $ h_t= h_x^{\frac{p+1}{2s-1}}$, where $p=2$ is the degree of the Lagrange polynomial in space. Preconditioners are approximated using one AMG V-cycle for each subsolve.}
  %\footnotesize
	\begin{tabular}{c| c| c| c c c c}
	\hline 
        %&     & GMRES   & GMRES      & GMRES        & GMRES \\
        %\multirow{6}{*}{s=2} &
       stage & $h_x^{-1}$  & DOF  & $\mathcal{P}_J$ & $\mathcal{P}_{GSL}$ & $\mathcal{P}_{DU}$ & $\mathcal{P}_{LD}$\Tstrut\Tstrut\Bstrut\Bstrut  \\
		%\cline{2-7}
		\hline
		%\hline
		\multirow{6}{*}{s=2}
		&8   & 450     & 14 (0.03) & 8 (0.01)   & 7 (0.01)  & 7 (0.01)\Tstrut  \\
		%\cline{2-7}
		&16  & 1922   & 14 (0.05)  & 8 (0.02)   & 7 (0.02)  & 7 (0.02)   \\
		%\cline{2-7}
		&32  & 7938    & 14 (0.10)  & 8 (0.06)   & 7 (0.06) & 7 (0.05)  \\
		%\cline{2-7}
		&64  & 32,258   & 14 (0.40) & 7 (0.19)   & 7 (0.19)  & 7 (0.19)   \\
		%\cline{2-7}
        &128 & 130,050  & 14 (1.48) & 7 (0.80)   & 7 (0.83)  & 7 (0.82)\Bstrut\\
    %\hline
	%\multicolumn{6}{|c|}{s=3}\\
	\hline 
	%\hline
	\multirow{5}{*}{s=3}
		&8   & 675     & 23 (0.06) & 10 (0.05)   & 10 (0.03)  & 9 (0.03)\Tstrut  \\
		%\cline{2-7}
		&16  & 2883   & 22 (0.06)  & 10 (0.03)   & 10 (0.03)  & 8 (0.02)  \\
		%\cline{2-7}
		&32  & 11,907  & 21 (0.17)  & 10 (0.09)   & 10 (0.08) & 8 (0.07)  \\
		%\cline{2-7}
		&64  & 48,387   & 22 (0.67) & 9 (0.29)   & 10 (0.32) & 8 (0.26)  \\
		%\cline{2-7}
        &128 & 195,075  & 21 (2.90) & 8 (1.12)   & 9 (1.25)  & 8 (1.12)\Bstrut \\
    \hline
    %\hline
    \multirow{5}{*}{s=4}
		&8   & 900     & 32 (0.10) & 13 (0.08)   & 13 (0.07)  & 10 (0.04)\Tstrut  \\
		%\cline{2-7}
		&16  & 3844    & 31 (0.11) & 12 (0.05)  & 13 (0.05)    & 10 (0.04) \\
		%\cline{2-7}
		&32  & 15,876   & 30 (0.31) & 12 (0.13)  & 12 (0.13)   & 10 (0.12) \\
		%\cline{2-7}
		&64  & 64,516   & 29 (1.16) & 11 (0.47)  & 12 (0.52) & 9 (0.39) \\
		%\cline{2-7}
        &128 & 260,100  & 27 (5.24)  & 11 (2.10) & 12 (2.29)  & 9 (1.73)\Bstrut\\ 
    \hline
     % \hline
    \multirow{5}{*}{s=5}
        &8   & 1125   & 42 (0.10) & 15 (0.08)  & 16 (0.08)  & 11 (0.06)\Tstrut \\
		%\cline{2-7}
		&16  & 4805   & 40 (0.23) & 15 (0.08)  & 16 (0.08)  & 11 (0.06) \\
		%\cline{2-7}
		&32  & 19,845 & 39 (0.50) & 13 (0.18)  & 15 (0.20)  & 11 (0.17) \\
		%\cline{2-7}
		&64  & 80,645 & 35 (1.80) & 13 (0.70)  & 15 (0.81)  & 11 (0.61) \\
		%\cline{2-7}
        &128 & 325,125 & 33 (8.44) & 12 (2.98) & 14 (3.49)  & 11 (2.75)\Bstrut\\ 
    \hline
    \multirow{5}{*}{s=6}
    &8 & 1350 & 53 (0.15) & 18 (0.08) & 19 (0.08) & 12 (0.04)\Tstrut \\
    &16 & 5766 & 50 (0.31) & 17 (0.11) & 18 (0.11) & 12 (0.08) \\
    &32 & 23,814 & 46 (0.73) & 16 (0.27) & 18 (0.29) & 12 (0.19) \\
    &64 &96,774 & 42 (2.66) & 15 (1.02) & 17 (1.14) & 12 (0.81) \\
    &128 & 390,150 & 38 (12.13) & 14 (4.27) & 17 (5.41) & 12 (3.56)\Bstrut \\
    \hline
    \multirow{5}{*}{s=7}
    &8 & 1575 & 62 (0.22) & 21 (0.10) & 23 (0.10) & 13 (0.09)\Tstrut \\
    &16 & 6727 & 57 (0.40) & 20 (0.15) & 22 (0.16) & 13 (0.09) \\
    &32 & 27,783 & 52 (0.98) & 19 (0.37) & 21 (0.40) & 13 (0.26) \\
    &64 & 112,903 & 48 (3.64) & 18 (1.49) & 20 (1.61) & 12 (0.98) \\
    &128 & 455,175 & 44 (17.06) & 17 (6.29) & 20 (7.24) & 12 (4.42)\Bstrut \\
    \hline
	\end{tabular}
	\label{tab:manufacturedResultsRadauIIA}
	\end{center}
	\end{table}

	% commented out old results for Lobatto IIIC with actual r.h.s 'b'
	\begin{comment}
	\begin{table}[htbp]
	    \centering
	      \caption{Iteration counts and elapsed time (times in seconds are shown in parentheses) for right-preconditioned GMRES to converge with relative residual tolerance $1.0 \times 10^{-8}$ for a $2D$ heat problem with various IRK methods with preconditioners $\mathcal{P}_{J}$, $\mathcal{P}_{GSL}$, $\mathcal{P}_{DU}$, and $\mathcal{P}_{LD}$. Here, $h_x=2^{-7}$ and $ h_t= h_x^{\frac{p+1}{q}}$, where $p=2$ is the degree of the Lagrange polynomial in space 
	      and $q$ is the order of the corresponding IRK method. 
	      Preconditioners are approximated using one AMG V-cycle for each subsolve.}
	    \begin{tabular}{c c c c c}
	    \hline
	        IRK method & $\mathcal{P}_{J}$ & $\mathcal{P}_{GSL}$ & $\mathcal{P}_{DU}$ & $\mathcal{P}_{LD}$\Tstrut\Tstrut\Bstrut\Bstrut \\
	        \hline
	         %G 2 & 12 (1.18) & 6 (0.58) & 7 (0.65) & 6 (0.54)\Tstrut\\
	         %G 3 & 18 (2.00) & 8 (1.09) & 10 (1.20) & 8 (0.95)\\
	         %G 4 & 25 (4.91) & 10 (1.66) & 13 (2.16) & 9 (1.49)\\
	         %G 5 & 32 (6.49) & 12 (2.43) & 16 (3.44) & 10 (2.03)\Bstrut\\
	         %\hline
	         L IIIC 2 & 12 (1.10) & 7 (0.67) & 6 (0.56) & 6 (0.59)\Tstrut\\
	         L IIIC 3 & 28 (3.30) & 12 (1.40) & 10 (1.13) & 10 (1.14)\\
	         L IIIC 4 & 41 (6.92) & 15 (2.34) & 15 (2.51) & 11 (1.81)\\
	         L IIIC 5 & 55 (12.87) & 18 (3.88) & 18 (3.94) & 12 (2.61)\Bstrut\\
	         \hline
	    \end{tabular}
	    \label{tab:resultsVariousIRK}
	\end{table}
	\end{comment}
	
We have also tested all four preconditioners for Lobatto IIIC methods. 
Table \ref{tab:manufacturedResultsLobattoIIIC} shows iteration counts and timing results for right-preconditioned GMRES to converge with a relative tolerance of $1.0 \times 10^{-8}$ for the $2D$ heat problem \eqref{eqn:testProblemHeat} using Lobatto IIIC methods in time and quadratic finite elements in space.
For the Lobatto IIIC methods, 
the $\mathcal{P}_{LD}$ preconditioner outperforms all of the other preconditioners both in iteration count and timing, with more significant improvements for larger problems and larger number of stages (higher order). 
%We did not consider Gauss--Legendre methods as Gauss--Legendre methods are not L-stable and are not suitable for the heat equation. 
%VEH - I think we say this earlier.
% Manufactured solution results for s-stage Lobatto IIIC
\begin{table}[htbp]
	\begin{center}
	 \caption{Iteration counts and elapsed time (times in seconds are shown in parentheses) for right-preconditioned GMRES to converge with relative residual tolerance $1.0 \times 10^{-8}$ for a $2D$ heat problem with $s$-stage Lobatto IIIC methods with preconditioners $\mathcal{P}_J$, $\mathcal{P}_{GSL}$, $\mathcal{P}_{DU}$, and $\mathcal{P}_{LD}$. Here we choose $ h_t= h_x^{\frac{p+1}{2s-2}}$, where $p=2$ is the degree of the Lagrange polynomial in space. Preconditioners are approximated using one AMG V-cycle for each subsolve.}
  %\footnotesize
	\begin{tabular}{c| c| c| c c c c}
	\hline 
       stage & $h_x^{-1}$  & DOF  & $\mathcal{P}_J$ & $\mathcal{P}_{GSL}$ & $\mathcal{P}_{DU}$ & $\mathcal{P}_{LD}$\Tstrut\Tstrut\Bstrut\Bstrut  \\
		\hline
		\multirow{6}{*}{s=2}
		&8   & 450     & 17 (0.01) & 10 (0.01)   & 7 (0.03)  & 7 (0.01)\Tstrut  \\
		%\cline{2-7}
		&16  & 1922   & 18 (0.03)  & 10 (0.02)   & 7 (0.02)  & 8 (0.02)   \\
		%\cline{2-7}
		&32  & 7938    & 19 (0.13)  & 10 (0.07)   & 8 (0.06) & 8 (0.06)  \\
		%\cline{2-7}
		&64  & 32,258   & 19 (0.47) & 10 (0.26)   & 7 (0.19)  & 8 (0.22)   \\
		%\cline{2-7}
        &128 & 130,050  & 18 (2.19) & 10 (1.21)   & 7 (0.86)  & 8 (1.00)\Bstrut\\
    %\hline
	%\multicolumn{6}{|c|}{s=3}\\
	\hline 
	%\hline
	\multirow{5}{*}{s=3}
		&8   & 675     & 27 (0.05) & 12 (0.03)   & 11 (0.02)  & 10 (0.02)\Tstrut  \\
		%\cline{2-7}
		&16  & 2883   & 29 (0.07)  & 12 (0.03)   & 11 (0.03)  & 10 (0.03)  \\
		%\cline{2-7}
		&32  & 11,907  & 27 (0.22)  & 12 (0.10)   & 11 (0.10) & 10 (0.08)  \\
		%\cline{2-7}
		&64  & 48,387   & 29 (0.96) & 11 (0.36)   & 11 (0.34) & 10 (0.35)  \\
		%\cline{2-7}
        &128 & 195,075  & 27 (3.89) & 11 (1.54)   & 10 (1.41)  & 9 (1.27)\Bstrut \\
    \hline
    %\hline
    \multirow{5}{*}{s=4}
		&8   & 900     & 40 (0.07) & 16 (0.05)   & 15 (0.05)  & 12 (0.05)\Tstrut  \\
		%\cline{2-7}
		&16  & 3844    & 40 (0.14) & 15 (0.06)  & 15 (0.06)    & 12 (0.05) \\
		%\cline{2-7}
		&32  & 15,876   & 39 (0.41) & 14 (0.15)  & 14 (0.15)   & 12 (0.13) \\
		%\cline{2-7}
		&64  & 64,516   & 38 (1.58) & 14 (0.60)  & 13 (0.57) & 11 (0.47) \\
		%\cline{2-7}
        &128 & 260,100  & 37 (7.58)  & 13 (2.46) & 13 (2.44)  & 11 (2.06)\Bstrut\\ 
    \hline
     % \hline
    \multirow{5}{*}{s=5}
        &8   & 1125   & 55 (0.08) & 19 (0.03)  & 19 (0.03)  & 13 (0.02)\Tstrut \\
		%\cline{2-7}
		&16  & 4805   & 55 (0.26) & 18 (0.09)  & 18 (0.08)  & 13 (0.07) \\
		%\cline{2-7}
		&32  & 19,845 & 51 (0.57) & 17 (0.20)  & 17 (0.22)  & 13 (0.16) \\
		%\cline{2-7}
		&64  & 80,645 & 49 (2.61) & 16 (0.85)  & 16 (0.84)  & 12 (0.64) \\
		%\cline{2-7}
        &128 & 325,125 & 45 (12.01) & 15 (3.64) & 15 (3.69)  & 12 (2.96)\Bstrut\\ 
    \hline
	\end{tabular}
	\label{tab:manufacturedResultsLobattoIIIC}
	\end{center}
\end{table}

Because a common time-stepping solution method for \eqref{eqn:testProblemHeat} is an adaptive strategy where 
the time step increases as the steady-state is approached, we also investigate the effectiveness of our 
preconditioner for a fixed spatial resolution while varying the time step size. 
In Table \ref{tab:varyTimeStepHeat}, we report iteration counts and timing results for left 
and right-preconditioned GMRES to converge with a relative residual tolerance of $1.0 \times 10^{-8}$ 
for the $2D$ heat problem with 2-stage and 7-stage Radau IIA methods. 
In this table, we have fixed the spatial resolution at $h_x^{-1} = 128$ and 
we vary the time step size $h_t$ between 0.05 and 5.0. Both the $\mathcal{P}_{GSL}$ preconditioner and our 
$\mathcal{P}_{LD}$ preconditioner are quite robust with respect to varying time step size.
Both preconditioners perform very similarly (and robustly) for the lower-order $s=2$ method; 
for the higher-order $s=7$ method, both preconditioners are robust with respect to the varying time step size, 
but our $\mathcal{P}_{LD}$  preconditioner outperforms $\mathcal{P}_{GSL}$ as a left preconditioner and as a right
preconditioner for all values of the time step.

% Varying time-step results for 2-D heat equation with Radau IIA
\begin{table}[htbp]
	\begin{center}
		\caption{Iteration counts and elapsed time (times in seconds are shown in parentheses) for left-preconditioned and right-preconditioned GMRES to converge with preconditioned relative residual tolerance $1.0 \times 10^{-8}$ for a $2D$ heat problem with $s$-stage Radau IIA methods with preconditioners  $\mathcal{P}_{GSL}$ and $\mathcal{P}_{LD}$. Here we keep $h_x^{-1}=128$ fixed and vary $h_t$ from 0.05 to 5.0. Preconditioners are approximated using one AMG V-cycle for each subsolve.}
		\begin{tabular}{c| c c c c c c}
			\hline
			\multicolumn{2}{c}{} &\multicolumn{2}{c}{left prec. GMRES} & & \multicolumn{2}{c}{right prec. GMRES}\Tstrut\Bstrut\\[0.5em]
			\cline{3-4}\cline{6-7}
			%\hline
			stage & $h_t$ & $\mathcal{P}_{GSL}$ & $\mathcal{P}_{LD}$ &&$\mathcal{P}_{GSL}$ &$\mathcal{P}_{LD}$ \Tstrut\Bstrut\\
			\hline
			\multirow{5}{*}{$s=2$}
			&0.05 & 8 (0.85) & 7 (0.71) & &7 (0.63) & 7 (0.63)  \Tstrut\\
			&0.1 & 8 (0.78) & 7 (0.73) & &6 (0.55) & 7 (0.63)\\
			&0.5 & 7 (0.74) & 7 (0.71) & &6 (0.56) & 7 (0.63)\\
			&1.0 & 7 (0.70) & 7 (0.70) & &6 (0.56) & 7 (0.64)\\
			&5.0 & 7 (0.72) & 8 (0.78) & &7 (0.63) & 7 (0.65)\Bstrut\\
			\hline
			\multirow{5}{*}{$s=7$}
			& 0.05 & 24 (7.38) & 14 (4.24) & &19 (5.43) & 12 (3.49)\Tstrut\\
			& 0.1 & 24 (7.76) & 14 (4.53) & &18 (5.62) & 12 (3.70)\\
			& 0.5 & 23 (7.42) & 14 (4.49) & &17 (5.34) & 13 (4.05)\\
			& 1.0 & 23 (6.77) & 14 (4.56) & &17 (4.89) & 13 (4.00)\\
			& 5.0 & 21 (6.74) & 15 (4.93) & &18 (5.64) &13 (4.17)\Bstrut\\
			\hline
		\end{tabular}
		\label{tab:varyTimeStepHeat}
	\end{center}
\end{table}

\subsubsection{Comparison results between $\widetilde{\mathcal{P}}_{GSL}$ with optimal coefficients and $\mathcal{P}_{LD}$ preconditioners}
	
 In this section we compare our $LDU$-based block lower triangular preconditioner $P_{LD}$ with the optimal block lower triangular Gauss--Seidel preconditioner introduced in  \cite{staff2006preconditioning}. 
 $\widetilde{\mathcal{P}}_{GSL}$ was developed using coefficients optimized to reduce the condition number of the left-preconditioned system $\widetilde{\mathcal{P}}_{GSL}^{-1} \mathcal{A}$.

In Table \ref{tab:comparisonRadauIIA}, we present results similar to those in Table \ref{tab:manufacturedResultsRadauIIA}. The problem set up is the same, but here we use left preconditioning for all of the preconditioners since the coefficients of $\widetilde{\mathcal{P}}_{GSL}$ were optimized for left preconditioning, and we include an additional column for $\widetilde{\mathcal{P}}_{GSL}$.
The first thing we note is that although the idea of optimizing coefficients to reduce the condition number of the preconditioned system is a perfectly sensible idea, for this problem the unoptimized $\mathcal{P}_{GSL}$ performs better than the optimized $\widetilde{\mathcal{P}}_{GSL}$ except on the smallest problems.  For all numbers of stages, $\mathcal{P}_{GSL}$ performs better than $\widetilde{\mathcal{P}}_{GSL}$ for problems sizes $h_x^{-1} = 32,64,$ and $128$.
We observe in this table that our $\mathcal{P}_{LD}$ has the best performance overall. Although it performs slightly worse as a left preconditioner than it did in Table \ref{tab:manufacturedResultsRadauIIA} as a right preconditioner, it nonetheless achieves the lowest iteration count and lowest timing of all of the preconditioners for all numbers of stages and all problem sizes.

 % Manufactured solution results for s-stage Radau IIA using left-preconditioned GMRES (comparison table)
 
\begin{table}[htbp]
	\begin{center}
	 \caption{Iteration counts and elapsed time (times in seconds are shown in parentheses) for left-preconditioned GMRES to converge with preconditioned relative residual tolerance $1.0 \times 10^{-8}$ for a $2D$ heat problem with $s$-stage Radau IIA methods with preconditioners $\mathcal{P}_J$, $\mathcal{P}_{GSL}$, $\widetilde{\mathcal{P}}_{GSL}$, $\mathcal{P}_{DU}$, and $\mathcal{P}_{LD}$. Here we choose $ h_t= h_x^{\frac{p+1}{2s-1}}$, where $p=2$ is the degree of the Lagrange polynomial in space. Preconditioners are approximated using one AMG V-cycle for each subsolve.}
  %\footnotesize
	\begin{tabular}{c| c| c| c c c c c}
	\hline 
        %&     & GMRES   & GMRES      & GMRES        & GMRES \\
        %\multirow{6}{*}{s=2} &
       stage & $h_x^{-1}$  & DOF  & $\mathcal{P}_J$ & $\mathcal{P}_{GSL}$  & $\widetilde{\mathcal{P}}_{GSL}$& $\mathcal{P}_{DU}$ & $\mathcal{P}_{LD}$\Tstrut\Tstrut\Bstrut\Bstrut  \\
		%\cline{2-7}
		\hline
		%\hline
		\multirow{6}{*}{s=2}
		&8   & 450     & 14 (0.02) & 9 (0.02) & 9 (0.01)  & 7 (0.01)  & 7 (0.01)\Tstrut  \\
		%\cline{2-7}
		&16  & 1922   & 16 (0.03)  & 9 (0.02) & 10 (0.02)   & 8 (0.02)  & 7 (0.02)   \\
		%\cline{2-7}
		&32  & 7938    & 16 (0.10)  & 9 (0.06) & 10 (0.06)  & 8 (0.05) & 7 (0.05)  \\
		%\cline{2-7}
		&64  & 32,258   & 17 (0.38) & 9 (0.24) & 10 (0.27)  & 8 (0.23)  & 7 (0.21)   \\
		%\cline{2-7}
        &128 & 130,050  & 18 (2.01) & 9 (1.16)  & 10 (1.24) & 8 (1.03)  & 7 (0.94)\Bstrut\\
    %\hline
	%\multicolumn{6}{|c|}{s=3}\\
	\hline 
	%\hline
	\multirow{5}{*}{s=3}
		&8   & 675     & 23 (0.06) & 11 (0.04) & 13 (0.03)  & 10 (0.04)  & 9 (0.03)\Tstrut  \\
		%\cline{2-7}
		&16  & 2883   & 25 (0.05)  & 12 (0.03) & 13 (0.03)  & 11 (0.03)  & 9 (0.02)  \\
		%\cline{2-7}
		&32  & 11,907  & 25 (0.20)  & 12 (0.09) & 14 (0.10)  & 11 (0.08) & 9 (0.07)  \\
		%\cline{2-7}
		&64  & 48,387   & 27 (0.77) & 12 (0.36) & 14 (0.43)  & 11 (0.34) & 9 (0.31)  \\
		%\cline{2-7}
        &128 & 195,075  & 28 (3.65) & 12 (1.65) & 14 (1.91)  & 11 (1.56)  & 9 (1.32)\Bstrut \\
    \hline
    %\hline
    \multirow{5}{*}{s=4}
		&8   & 900     & 32 (0.07) & 15 (0.06) & 18 (0.05)  & 14 (0.06)  & 11 (0.04)\Tstrut  \\
		%\cline{2-7}
		&16  & 3844    & 34 (0.13) & 15 (0.06) & 19 (0.07) & 14 (0.06)    & 11 (0.05) \\
		%\cline{2-7}
		&32  & 15,876   & 36 (0.38) & 15 (0.17) &19 (0.21) & 14 (0.16)   & 11 (0.13) \\
		%\cline{2-7}
		&64  & 64,516   & 36 (1.52) & 15 (0.70) & 19 (0.84) & 14 (0.64) & 11 (0.51) \\
		%\cline{2-7}
        &128 & 260,100  & 36 (6.51)  & 15 (2.88) & 19 (3.63) & 14 (2.74)  & 11 (2.20)\Bstrut\\ 
    \hline
     % \hline
    \multirow{5}{*}{s=5}
        &8   & 1125   & 44 (0.11) & 18 (0.08) & 22 (0.06)  & 17 (0.07)  & 12 (0.06)\Tstrut \\
		%\cline{2-7}
		&16  & 4805   & 45 (0.21) & 18 (0.08) & 22 (0.10) & 17 (0.08)  & 12 (0.06) \\
		%\cline{2-7}
		&32  & 19,845 & 46 (0.61) & 18 (0.26) & 22 (0.28)  & 17 (0.23)  & 12 (0.17) \\
		%\cline{2-7}
		&64  & 80,645 & 47 (2.21) & 18 (0.90) & 22 (1.18)  & 17 (0.93)  & 12 (0.74) \\
		%\cline{2-7}
        &128 & 325,125 & 47 (11.15) & 18 (4.41) & 22 (5.45) & 17 (4.26)  & 12 (3.00)\Bstrut\\ 
    \hline
    \multirow{5}{*}{s=6}
    &8 & 1350 & 55 (0.14) & 22 (0.08) & 24 (0.07) & 20 (0.07) & 13 (0.04)\Tstrut \\
    &16 & 5766 & 55 (0.32) & 22 (0.14) & 24 (0.16) & 21 (0.14) & 14 (0.09) \\
    &32 & 23,814 & 56 (0.93) & 22 (0.38) & 24 (0.40) & 20 (0.34) & 14 (0.24) \\
    &64 &96,774 & 56 (3.77) & 21 (1.49) &24 (1.69) & 20 (1.42) & 13 (0.94) \\
    &128 & 390,150 & 56 (15.89) & 21 (6.21) &23 (6.89) & 20 (5.63) & 13 (3.98)\Bstrut \\
    \hline
	\end{tabular}
	\label{tab:comparisonRadauIIA}
	\end{center}
\end{table}
   
In Table \ref{tab:comparisonLobattoIIIC}, we present results 
similar to those in Table \ref{tab:comparisonRadauIIA} but using Lobatto IIIC methods. The problem set up is again the same, using left preconditioning for all of the preconditioners. 
The results are very similar to those in the previous table.
The unoptimized $\mathcal{P}_{GSL}$ performs better than the optimized $\widetilde{\mathcal{P}}_{GSL}$ on all of these problems. 
Our $\mathcal{P}_{LD}$ preconditioner has the best performance overall.

% Manufactured solution results for s-stage Lobatto IIIC using left-preconditioned GMRES (comparison table)
\begin{table}[htbp]
	\begin{center}
	 \caption{Iteration counts and elapsed time (times in seconds are shown in parentheses) for left-preconditioned GMRES to converge with preconditioned relative residual tolerance $1.0 \times 10^{-8}$ for a $2D$ heat problem with $s$-stage Lobatto IIIC methods with preconditioners $\mathcal{P}_J$, $\mathcal{P}_{GSL}$, $\widetilde{\mathcal{P}}_{GSL}$, $\mathcal{P}_{DU}$, and $\mathcal{P}_{LD}$. Here we choose $ h_t= h_x^{\frac{p+1}{2s-2}}$, where $p=2$ is the degree of the Lagrange polynomial in space. Preconditioners are approximated using one AMG V-cycle for each subsolve.}
  %\footnotesize
	\begin{tabular}{c| c| c| c c c c c}
	\hline 
       stage & $h_x^{-1}$  & DOF  & $\mathcal{P}_J$ & $\mathcal{P}_{GSL}$  & $\widetilde{\mathcal{P}}_{GSL}$& $\mathcal{P}_{DU}$ & $\mathcal{P}_{LD}$\Tstrut\Tstrut\Bstrut\Bstrut  \\
		\hline
		\multirow{6}{*}{s=2}
		&8   & 450     & 18 (0.02) & 11 (0.02) &18 (0.02)   & 7 (0.02)  & 8 (0.02)\Tstrut  \\
		&16  & 1922   & 20 (0.04)  & 12 (0.02) &20 (0.03)  & 8 (0.02)  & 8 (0.02)   \\
		&32  & 7938    & 21 (0.15)  & 12 (0.09) &21 (0.14)  & 8 (0.06) & 8 (0.06)  \\
		&64  & 32,258   & 21 (0.56) & 12 (0.34) &21 (0.56)  & 8 (0.24) & 8 (0.23)   \\
        &128 & 130,050  & 21 (2.40) & 12 (1.49) &21 (2.50) & 8 (1.05) & 8 (1.05)\Bstrut\\
    %\hline
	%\multicolumn{6}{|c|}{s=3}\\
	\hline 
	%\hline
	\multirow{5}{*}{s=3}
		&8   & 675     & 27 (0.06) & 14 (0.04) &25 (0.04) & 11 (0.03)  & 11 (0.03)\Tstrut  \\
		&16  & 2883   & 32 (0.08)  & 14 (0.04) &27 (0.06)  & 12 (0.03)  & 11 (0.03)  \\
		&32  & 11,907  & 33 (0.23)  & 15 (0.11) &28 (0.20)  & 12 (0.09) & 11 (0.09)  \\
		&64  & 48,387   & 35 (1.06) & 15 (0.49) &30 (0.98)  & 12 (0.40) & 11 (0.37)  \\
        &128 & 195,075  & 37 (3.89) & 15 (1.87) &30 (3.76)  & 12 (1.67)  & 11 (1.58)\Bstrut \\
    \hline
    %\hline
    \multirow{5}{*}{s=4}
		&8   & 900     & 41 (0.09) & 18 (0.06) &31 (0.07)  & 15 (0.06)  & 13 (0.05)\Tstrut  \\
		&16  & 3844    & 44 (0.13) & 18 (0.06) &35 (0.11) & 16 (0.06)    & 13 (0.05) \\
		&32  & 15,876   & 46 (0.51) & 18 (0.20) &35 (0.38)  & 16 (0.18)   & 13 (0.15) \\
		&64  & 64,516   & 49 (2.14) & 18 (0.80) &37 (1.66) & 16 (0.74) & 13 (0.60) \\
        &128 & 260,100  & 49 (8.93)  & 18 (3.23) &36 (6.60) & 16 (3.11)  & 13 (2.53)\Bstrut\\ 
    \hline
	\end{tabular}
	\label{tab:comparisonLobattoIIIC}
	\end{center}
\end{table}

It should be noted that since the results in Tables \ref{tab:comparisonRadauIIA} and \ref{tab:comparisonLobattoIIIC} are using left preconditioning, the relative residuals are actually preconditioned relative residuals. So the norms of the various preconditioners could be having an effect on the iteration counts. Since all of our numerical tests were run using the method of manufactured solutions, in Tables \ref{tab:errorLeftPrec} and \ref{tab:errorRightPrec} we examine the relative errors for each method. 

In Table \ref{tab:errorLeftPrec} we report the relative errors for various numbers of stages for Radau IIA and Lobatto IIIC with the five preconditioners $\mathcal{P}_J$, $\mathcal{P}_{GSL}$,  $\widetilde{\mathcal{P}}_{GSL}$, $\mathcal{P}_{DU}$, and  $\mathcal{P}_{LD}$ all used as left preconditioners. The problem size is fixed at $h_x^{-1} = 128$, which is the largest of the problems we considered, and $ h_t= h_x^{\frac{p+1}{q}}$, where $p=2$ is the degree of the Lagrange polynomial in space and $q$ is the order of the corresponding IRK method
(that is, $q = 2s-1$ for Radau IIA and $q = 2s -2$ for Lobatto IIIC).
The stopping criteria is the same as in all previous tables. 
With left preconditioning, problems with Radau IIA and Lobatto IIIC methods and all numbers of stages yielded very similar relative errors.

%$\mathcal{P}_J$ 
%$\mathcal{P}_{GSL}$ 
%$\widetilde{\mathcal{P}}_{GSL}$ 
%$\mathcal{P}_{DU}$ 
%$\mathcal{P}_{LD}$ 
  
\begin{table}[htbp]
       \centering
       \caption{Relative error for left-preconditioned GMRES converging to a preconditioned relative residual tolerance of $1.0 \times 10^{-8}$ for a $2D$ heat problem with various IRK methods with preconditioners $\mathcal{P}_J$, $\mathcal{P}_{GSL}$, $\widetilde{\mathcal{P}}_{GSL}$, $\mathcal{P}_{DU}$, and $\mathcal{P}_{LD}$. Here, $h_x=2^{-7}$ and $ h_t= h_x^{\frac{p+1}{q}}$, where $p=2$ is the degree of the Lagrange polynomial in space and $q$ is the order of the corresponding IRK method (that is, $q=2s-1$ for Radau IIA and $q=2s-2$ for Lobatto IIIC). Preconditioners are approximated using one AMG V-cycle for each subsolve.}
       \begin{tabular}{c c c c c c}
       \hline
       IRK method & $\mathcal{P}_J$ & $\mathcal{P}_{GSL}$  & $\widetilde{\mathcal{P}}_{GSL}$& $\mathcal{P}_{DU}$ & $\mathcal{P}_{LD}$\Tstrut\Tstrut\Bstrut\Bstrut  \\
       \hline
        R IIA 2 & $7.1\times 10^{-9}$ & $6.1\times 10^{-9}$ & $1.4\times 10^{-9}$ & $1.8\times 10^{-9}$ & $4.0\times 10^{-9}$\Tstrut \\
        R IIA 3 & $8.5\times 10^{-9}$ & $6.9\times 10^{-9}$ & $7.9\times 10^{-9}$ & $9.4\times 10^{-9}$ & $6.4\times 10^{-9}$\\
        R IIA 4 & $4.6\times 10^{-8}$ & $8.5\times 10^{-9}$ & $9.1\times 10^{-9}$ & $3.8\times 10^{-8}$ & $2.9\times 10^{-9}$ \\
        R IIA 5 & $5.1\times 10^{-8}$ & $9.0\times 10^{-9}$ & $8.1\times 10^{-9}$ & $4.6\times 10^{-8}$ & $5.4\times 10^{-9}$ \\
        R IIA 6 & $7.7\times 10^{-8}$ & $1.0\times 10^{-8}$ & $8.3\times 10^{-9}$ & $9.8\times 10^{-8}$ & $7.9\times 10^{-9}$\Bstrut \\
        \hline
        L IIIC 2 & $8.0\times 10^{-9}$ & $1.8\times 10^{-9}$ & $8.0\times 10^{-9}$ & $3.9\times 10^{-9}$ & $6.6\times 10^{-9}$\Tstrut \\
        L IIIC 3 & $2.0\times 10^{-8}$ & $6.7\times 10^{-9}$ & $1.0\times 10^{-8}$ & $1.8\times 10^{-8}$ & $1.1\times 10^{-8}$\\
        L IIIC 4 & $5.7\times 10^{-8}$ & $1.2\times 10^{-8}$ & $7.9\times 10^{-9}$ & $2.3\times 10^{-8}$ & $8.9\times 10^{-9}$\Bstrut \\
        \hline
       \end{tabular}
       \label{tab:errorLeftPrec}
   \end{table}
  
In Table \ref{tab:errorRightPrec} we report the relative errors for various numbers of stages for Radau IIA and Lobatto IIIC with the four preconditioners $\mathcal{P}_J$, $\mathcal{P}_{GSL}$,  $\mathcal{P}_{DU}$, and  $\mathcal{P}_{LD}$ this time all used as right preconditioners. We do not include the optimized GSL-based preconditioner $\widetilde{\mathcal{P}}_{GSL}$ since it is optimized for left preconditioning. The problem size is again fixed at $h_x^{-1} = 128$, and as before $ h_t= h_x^{\frac{p+1}{q}}$, where $p=2$ is the degree of the Lagrange polynomial in space and $q$ is the order of the corresponding IRK method (that is, $q=2s-1$ for Radau IIA and $q = 2s-2$ for Lobatto IIIC). The stopping criteria is the same as in all previous tables. 
With right reconditioning, our $\mathcal{P}_{LD}$ preconditioner achieves almost an order of magnitude better relative error on almost all of the runs.
  
\begin{table}[htbp]
       \centering
       \caption{Relative error for right-preconditioned GMRES converging to a relative residual tolerance of $1.0 \times 10^{-8}$ for a $2D$ heat problem with various IRK methods with preconditioners $\mathcal{P}_J$, $\mathcal{P}_{GSL}$, $\mathcal{P}_{DU}$, and $\mathcal{P}_{LD}$. Here, $h_x=2^{-7}$ and $ h_t= h_x^{\frac{p+1}{q}}$, where $p=2$ is the degree of the Lagrange polynomial in space and $q$ is the order of the corresponding IRK method (that is, $q=2s-1$ for Radau IIA and $q=2s-2$ for Lobatto IIIC). Preconditioners are approximated using one AMG V-cycle for each subsolve.}
       \begin{tabular}{c c c c c c}
       \hline
       IRK method & $\mathcal{P}_J$ & $\mathcal{P}_{GSL}$  & $\mathcal{P}_{DU}$ & $\mathcal{P}_{LD}$\Tstrut\Tstrut\Bstrut\Bstrut  \\
       \hline
        R IIA 2 & $1.5\times 10^{-6}$ & $1.4\times 10^{-6}$ & $2.8\times 10^{-7}$ & $3.1\times 10^{-8}$\Tstrut \\
        R IIA 3 & $1.5\times 10^{-5}$ & $1.5\times 10^{-5}$ & $9.0\times 10^{-6}$ & $2.8\times 10^{-7}$\\
        R IIA 4 & $2.7\times 10^{-5}$ & $4.6\times 10^{-6}$ & $8.3\times 10^{-6}$ & $6.9\times 10^{-7}$\\
        R IIA 5 & $3.8\times 10^{-5}$ & $1.1\times 10^{-5}$ & $1.5\times 10^{-5}$ & $1.3\times 10^{-7}$ \\
        R IIA 6 & $5.6\times 10^{-5}$ & $2.0\times 10^{-5}$ & $9.5\times 10^{-6}$ & $8.2\times 10^{-8}$\\
        R IIA 7 & $1.7\times 10^{-6}$ & $3.9\times 10^{-7}$ & $8.3\times 10^{-7}$ & $2.4\times 10^{-7}$\Bstrut \\
        \hline
        L IIIC 2 & $2.6\times 10^{-7}$ & $1.3\times 10^{-7}$ & $2.5\times 10^{-7}$ & $2.1\times 10^{-8}$\Tstrut \\
        L IIIC 3 & $1.2\times 10^{-5}$ & $8.1\times 10^{-6}$ & $9.2\times 10^{-6}$ & $6.9\times 10^{-6}$\\
        L IIIC 4 & $2.5\times 10^{-5}$ & $1.3\times 10^{-5}$ & $1.6\times 10^{-5}$ & $1.4\times 10^{-6}$\\
        L IIIC 5 & $6.1\times 10^{-5}$ & $1.3\times 10^{-5}$ & $4.0\times 10^{-5}$ & $1.2\times 10^{-6}$\Bstrut \\
        \hline
       \end{tabular}
       \label{tab:errorRightPrec}
\end{table}
  
\subsection{Advection-Diffusion Equation Results}
\label{sec:NumericalAD}
In this section, we examine the performance of the preconditioners 
$\mathcal{P}_{GSL}$ and $\mathcal{P}_{LD}$, on the $2D$ double-glazing advection-diffusion problem 
\eqref{eqn:testProblemAD}. For the double-glazing problem, we only present results for the
$\mathcal{P}_{GSL}$ and $\mathcal{P}_{LD}$ preconditioners since in all of our previous experiments 
they performed best. We apply both as left preconditioners in these results, but our results are
similar when they are applied as right preconditioners. 
%\textcolor{red}{This needs to be written}

In Table \ref{tab:resultsDG}, we report iteration counts and elapsed time for left-preconditioned 
GMRES to converge with preconditioned relative residual tolerance $1.0 \times 10^{-8}$ for 
the double-glazing problem \eqref{eqn:testProblemAD}
with $s$-stage Radau IIA methods using the preconditioners $\mathcal{P}_{GSL}$ and $\mathcal{P}_{LD}$. 
We vary the number of stages from $s=2$ to $s=7$ for problems with spatial resolution $h_x^{-1} = 64$ and $128$.
Here we choose $ h_t= h_x^{\frac{p+1}{2s-1}}$, where $p=1$ is the degree of the Lagrange polynomial in space.
We present results for $\epsilon = 0.04$, which is diffusion dominated and $\epsilon = 0.005$, which is weakly 
advection dominated. 
Preconditioners are approximated using one AMG V-cycle for each subsolve.

For the lowest-order problems with $s=2$, $\mathcal{P}_{GSL}$ and our $\mathcal{P}_{LD}$ preconditioner
perform similarly. For the diffusion dominated problem, $\epsilon = 0.04$, $\mathcal{P}_{LD}$ outperforms 
$\mathcal{P}_{GSL}$ for all of the higher-order problems, $s=3$ through $s=7$, with the improved performance
being more pronounced on the highest-order problems. For the weakly advection-dominated problem, 
$\epsilon = 0.005$, the results are more mixed. For the smaller problems of each order, $h_x^{-1} = 64$,
$\mathcal{P}_{LD}$ takes fewer iterations in most cases, but takes longer in compute time in a few cases.
For the largest higher-order problems, $s=4$ through $s=7$ with $h_x^{-1} = 128$, 
$\mathcal{P}_{LD}$ performs better than $\mathcal{P}_{GSL}$ in both iteration count and timing.

%============================================================================
%% Results for the Double-Glazing problem
%============================================================================

\begin{table}[htbp]
	\begin{center}
			 \caption{Iteration counts and elapsed time (times in seconds are shown in parentheses) for left-preconditioned GMRES to converge with preconditioned relative residual tolerance $1.0 \times 10^{-8}$ for the double-glazing problem with $s$-stage Radau IIA methods with preconditioners  $\mathcal{P}_{GSL}$ and $\mathcal{P}_{LD}$. Here we choose $ h_t= h_x^{\frac{p+1}{2s-1}}$, where $p=1$ is the degree of the Lagrange polynomial in space. Preconditioners are approximated using one AMG V-cycle for each subsolve.}
		    \begin{tabular}{c| c| c| c| c c| c c}
			%\cline{5-8}
			\hline
			 \multicolumn{4}{c|}{} &\multicolumn{2}{c|}{$\epsilon = 0.04$} &  \multicolumn{2}{c}{$\epsilon = 0.005$}\Tstrut\Bstrut\\
			%\cline{4-7}
			\hline
			stage & $h_x^{-1}$ & $h_t$& $DOF$ & $\mathcal{P}_{GSL}$ & $\mathcal{P}_{LD}$ &$\mathcal{P}_{GSL}$ &$\mathcal{P}_{LD}$ \Tstrut\Bstrut\\
			\hline
			\multirow{2}{*}{$s=2$}
			&64 & 0.0625 & 33282 & 9 (0.50) & 7 (0.43) & 11 (0.89) & 11 (0.83) \Tstrut\\
			&128 & 0.0394 & 132098 & 9 (1.72) & 7 (1.44) & 10 (3.74) & 10 (3.84)\Bstrut\\
			\hline
			\multirow{2}{*}{$s=3$}
			&64 & 0.1895 & 49923 & 14 (1.13) & 10 (0.84) & 17 (1.92) & 16 (1.83) \Tstrut\\
			&128 & 0.1436 & 198147 & 14 (3.96) & 9 (2.73) & 15 (8.23) & 14 (9.71)\Bstrut\\
			\hline
			\multirow{2}{*}{$s=4$}
			&64 &0.3048 &66564 & 18 (1.95) & 12 (1.28) & 22 (3.22) & 21 (3.04)  \Tstrut\\
			&128 &0.25 & 264196 & 18 (6.47) & 11 (4.44) & 20 (16.83) & 17 (15.17)\Bstrut\\
			\hline
			\multirow{2}{*}{$s=5$}
			&64 & 0.3969 & 83205 & 22 (2.62) & 13 (1.58) & 26 (4.54) & 25 (4.29) \Tstrut\\
			&128 & 0.3402 & 330245 & 22 (10.38) & 13 (6.45) & 25 (22.75) & 21 (22.33)\Bstrut\\
			\hline
			\multirow{2}{*}{$s=6$}
			&64 & 0.4695 & 99846 & 25 (3.94) & 14 (2.33) & 31 (6.34) & 28 (6.51) \Tstrut\\
			&128 & 0.4139 & 396294 & 27 (15.84) & 14 (10.59) & 30 (34.10) & 23 (33.38)\Bstrut \\
			\hline
			\multirow{2}{*}{$s=7$}
			&64 & 0.5274 & 116487 & 30 (5.39) & 16 (2.81) & 36 (8.53) & 30 (9.04) \Tstrut\\
			&128 & 0.4740 & 462343 & 31 (20.89) & 15 (10.34) & 36 (55.84) & 25 (36.33)\Bstrut\\
			\hline
		\end{tabular}
		\label{tab:resultsDG}
	\end{center}
\end{table}

As we did for the heat equation, in Table \ref{tab:varyTimeStepDG}, 
we investigate the robustness of the preconditioners for a fixed spatial resolution while varying
the time step size. In this table, we present iteration counts and timing results for left-preconditioned
GMRES to converge with a preconditioned relative residual tolerance of 
$1.0 \times 10^{-8}$ for the $2D$ double-glazing problem with $s$-stage Radau IIA methods using 
the preconditioners $\mathcal{P}_{GSL}$ and $\mathcal{P}_{LD}$. We present results for a lower-order method $s=2$ 
and a higher-order method $s=7$ and for  $\epsilon = 0.04$, which is diffusion dominated, 
and $\epsilon = 0.005$, which is weakly  advection dominated. 
In this table, we keep $h_x^{-1}=128$ fixed and vary the time step size $h_t$ from 0.05 to 5.0. 
Preconditioners are approximated using one AMG V-cycle for each subsolve.

For the diffusion-dominated problem with the lower-order 
method, our $\mathcal{P}_{LD}$ preconditioner mildly outperforms $\mathcal{P}_{GSL}$ for all time step sizes,
and both preconditioners are robust with respect to varying the time step size. 
For the higher-order method, both preconditioners take more iterations and more time as the step 
size increases, but $\mathcal{P}_{LD}$ grows less dramatically than $\mathcal{P}_{GSL}$ and takes fewer
iterations and less time for all time step sizes.

For the weakly advection-dominated problem, $\epsilon = 0.005$, $\mathcal{P}_{LD}$ 
performs very similarly to $\mathcal{P}_{GSL}$ on the lower-order problem, $s=2$, 
and both preconditioners are robust with respect to the varying time step size.
For the higher-order method, both preconditioners take more iterations and more time as the step 
size increases, but $\mathcal{P}_{LD}$ grows much less dramatically than $\mathcal{P}_{GSL}$ and takes fewer
iterations and less time for all time step sizes.

\begin{table}[htbp]
	\begin{center}
		\caption{Iteration counts and elapsed time (times in seconds are shown in parentheses) for left-preconditioned GMRES to converge with preconditioned relative residual tolerance $1.0 \times 10^{-8}$ for a $2D$ double-glazing problem with $s$-stage Radau IIA methods with preconditioners  $\mathcal{P}_{GSL}$ and $\mathcal{P}_{LD}$. Here we keep $h_x^{-1}=128$ fixed and vary $h_t$ from 0.05 to 5.0. Preconditioners are approximated using one AMG V-cycle for each subsolve.}
		\begin{tabular}{c| c c c c c c}
			\hline
			\multicolumn{2}{c}{} &\multicolumn{2}{c}{$\epsilon = 0.04$} & & \multicolumn{2}{c}{$\epsilon = 0.005$}\Tstrut\Bstrut\\[0.5em]
			\cline{3-4}\cline{6-7}
			%\hline
			stage & $h_t$ & $\mathcal{P}_{GSL}$ & $\mathcal{P}_{LD}$ &&$\mathcal{P}_{GSL}$ &$\mathcal{P}_{LD}$ \Tstrut\Bstrut\\
			\hline
			\multirow{5}{*}{$s=2$}
			&0.05 & 9 (1.85) & 7 (1.50) & &11 (3.95) & 10 (4.12)  \Tstrut\\
			&0.1 & 10 (1.96) & 7 (1.43) & &11 (4.24) & 11 (5.01)\\
			&0.5 & 11 (2.05) & 8 (1.76) & &12 (4.64) & 11 (4.64)\\
			&1.0 & 10 (2.08) & 7 (1.47) & &13 (4.76) & 11 (4.72)\\
			&5.0 & 10 (1.95) & 8 (1.69) & &12 (5.63) & 12 (5.43)\Bstrut\\
			\hline
			\multirow{5}{*}{$s=7$}
			& 0.05 & 26 (18.39) & 13 (9.78) & &27 (41.80) & 18 (28.56)\Tstrut\\
			& 0.1 & 27 (22.09) & 14 (9.80) & &28 (43.18) & 21 (29.72)\\
			& 0.5 & 32 (24.97) & 15 (12.29) & &36 (53.60) & 25 (40.76)\\
			& 1.0 & 37 (30.81) & 17 (13.88) & &45 (76.24) & 26 (36.33)\\
			& 5.0 & 46 (35.33) & 20 (16.14) &\, &68 (109.26) & 28 (36.66)\Bstrut\\
			\hline
		\end{tabular}
		\label{tab:varyTimeStepDG}
	\end{center}
\end{table}
% Removed NS section. I assume we will save that work for a follow-up paper.

%\clearpage 
\section{Conclusion}
%Explicit time integrators for parabolic PDE
%are subject to the restrictive timestep limit $h_t \lesssim h_x^2$, so A-stable integrators are essential. 
%based on an LDU factorization. Solves on individual blocks are accomplished using a multigrid algorithm. 
In this paper, we have introduced a new preconditioner 
for the large, structured systems appearing in 
implicit Runge--Kutta time integration of parabolic PDE problems. 
Our preconditioner is based on a block $LDU$ factorization, and
for scalability, we have used a single AMG V-cycle on all subsolves. 
We have compared our preconditioner $\mathcal{P}_{LD}$ in condition number and eigenvalue distribution to other preconditioners, 
and have 
demonstrated its effectiveness on two test problems, the heat equation and the double-glazing advection-diffusion problem. We have found that it is
scalable ($h_x$-independent) and yields better timing results than other preconditioners currently in the literature:
block Jacobi, block Gauss--Seidel, and the optimized block Gauss--Seidel method of~\cite{staff2006preconditioning}. It is also robust with respect to varying time step sizes for a fixed spatial resolution. We ran experiments with implicit Runge--Kutta stages up to $s=7$, and have found that the new preconditioner outperforms the others, with the improvement becoming more pronounced as spatial discretization is refined and as temporal order is increased. 
%\begin{itemize}
%\item With preconditioning (which seems promising), we can use IRK
      %methods and have higher order with L-stability. 
%\item Leverage existing preconditioners (PCD, LSC, etc.) on the subblocks.
%\item We have more methods to try (especially for $s>2$).
%\item Combine these ideas with W-transform techniques of L. O. Jay \& T.
      %Braconnier (1999).  (Block LU of a tridiagonal system. $P = \hat{L} \hat{U}$.
%\end{itemize}

%\newpage
\bibliographystyle{siamplain}
\bibliography{refs}

\begin{thebibliography}{1}

\bibitem{elman2005fem}
{\sc H.~Elman, D.~Silvester, and A.~Wathen}, {\em Finite Elements and Fast
  Iterative Solvers}, Oxford University Press, 2005.

\bibitem{hairer1993solving}
{\sc E.~Hairer, S.~P. N{\o}rsett, and G.~Wanner}, {\em Solving Ordinary
  Differential Equations I, Nonstiff Problems}, Springer Berlin Heidelberg,
  1993.

\bibitem{juntunen2009nitsche}
{\sc M.~Juntunen and R.~Stenberg}, {\em Nitsche’s method for general boundary
  conditions}, Mathematics of computation, 78 (2009), pp.~1353--1374.

\bibitem{long_unified_2010}
{\sc K.~Long, R.~Kirby, and B.~van Bloemen~Waanders}, {\em Unified embedded
  parallel finite element computations via software-based {Fr{\'e}chet}
  differentiation}, SIAM Journal on Scientific Computing, 32 (2010),
  pp.~3323--3351, \url{https://doi.org/10.1137/09076920X},
  \url{http://epubs.siam.org/doi/10.1137/09076920X} (accessed 2019-07-29).

\bibitem{mardal2007order}
{\sc K.-A. Mardal, T.~K. Nilssen, and G.~A. Staff}, {\em Order-optimal
  preconditioners for implicit {Runge}--{Kutta} schemes applied to parabolic
  {PDEs}}, SIAM Journal on Scientific Computing, 29 (2007), pp.~361--375.

\bibitem{ifiss}
{\sc D.~Silvester, H.~Elman, and A.~Ramage}, {\em {I}ncompressible {F}low and
  {I}terative {S}olver {S}oftware ({IFISS}) version 3.5}, September 2016.
\newblock {\tt http://www.manchester.ac.uk/ifiss/}.

\bibitem{staff2006preconditioning}
{\sc G.~A. Staff, K.-A. Mardal, and T.~K. Nilssen}, {\em Preconditioning of
  fully implicit {Runge}--{Kutta} schemes for parabolic {PDEs}}, Modeling,
  Identification and Control, 27 (2006), pp.~109--123,
  \url{https://doi.org/10.4173/mic.2006.2.3}.

\bibitem{stenberg1995some}
{\sc R.~Stenberg}, {\em On some techniques for approximating boundary
  conditions in the finite element method}, Journal of Computational and
  applied Mathematics, 63 (1995), pp.~139--148.

\bibitem{wanner1996solving}
{\sc G.~Wanner and E.~Hairer}, {\em Solving Ordinary Differential Equations II:
  Stiff Problems}, Springer Berlin Heidelberg, 1996.

\end{thebibliography}

\end{document}